\documentclass[11pt]{article} 

\usepackage{amsfonts}
\usepackage{amscd}
\usepackage{amssymb}
\usepackage{amsthm}
\usepackage{amsmath}
\usepackage{comment}
\usepackage{stmaryrd}
\usepackage{amsfonts}
\usepackage{mathrsfs}
\usepackage{color}
\usepackage{youngtab}
\usepackage [all]{xy}
\usepackage{rotating}

\theoremstyle{plain}
\newtheorem{thm}{Theorem}[section]

\theoremstyle{definition}

\theoremstyle{example}

\theoremstyle{remark}

\numberwithin{equation}{section}

\theoremstyle{plain}

\renewcommand{\ge}{\geqslant}

\renewcommand{\le}{\leqslant}

\def\cL{\mathcal{L}}

\def\AA{\mathbb{A}}

\def\CC{\mathbb{C}}

\def\LL{\mathbb{L}}

\def\RR{\mathbb{R}}

\def\ZZ{\mathbb{Z}}

\def\fh{\mathfrak{h}}

\def\Hom{\mathrm{Hom}}

\def\im{\mathrm{im}}
\def\Ind{\mathrm{Ind}}

\def\pt{\mathrm{pt}}

\def\mapright#1{\smash{\mathop
        {\longrightarrow}\limits^{#1}}}

\def\mapup#1{\Big\uparrow
   \rlap{$\vcenter{\hbox{$\scriptstyle#1$}}$}}

\makeatletter
\renewcommand{\@makefnmark}{\mbox{\textsuperscript{}}}
\makeatother

\title{Generalized Schubert Calculus}
\author{
Nora Ganter\ and\ Arun Ram \\
Department of Mathematics and Statistics \\
University of Melbourne \\
Parkville VIC 3010 Australia \\
nganter@unimelb.edu.au \ and\ 
aram@unimelb.edu.au \\
\\
{\sl Dedicated to C.S. Seshadri on the occasion of his 80th birthday}
}
\date{}

\begin{document}

\maketitle


\begin{abstract}
\noindent
In this paper we study the $T$-equivariant generalized cohomology of flag varieties using two models,
the Borel model and the moment graph model.  We study the differences between the Schubert classes and the 
Bott-Samelson classes.  After setup of the general framework we compute, for classes of Schubert varieties of
complex dimension $\le 3$ in rank 2 (including $A_2$, $B_2$, $G_2$ and $A_1^{(1)}$), moment graph representatives, 
Pieri-Chevalley formulas and products of Schubert classes.  These computations generalize the computations in 
equivariant K-theory for rank 2 cases which are given in Griffeth-Ram \cite{GR}.
\end{abstract}

\footnote{AMS Subject Classifications: Primary 14M17; Secondary 14N15.}

\section{Introduction}

This paper is a study of the generalized equivariant cohomology of flag varieties.
We set up a general framework for working with the generalized (equivariant) Schubert calculus which allows for detailed
study without the need for knowledge of cobordism or generalized cohomology theories.
Working in the context of a complex reductive algebraic group $G$, the 
(generalized) flag variety is $G/B$, where $B$ is a Borel subgroup containing the maximal torus $T$.  
The equivariant generalized cohomology theory $h_T$ comes with a (formal) group which is
used to combinatorially construct the ring $S=h_T(\pt)$.  The Borel model presents $h_T(G/B)$
as a `coinvariant ring' $S\otimes_{S^{W_0}} S$ and the moment graph model presents $h_T(G/B)$ via the image
of the inclusions of the $T$-fixed points of $G/B$.  Special cases of generalized equivariant cohomology theories
are `ordinary' cohomology (corresponding to the additive group) and $K$-theory (corresponding to the
multiplicative group).  The universal formal group law corresponds to complex cobordism.


Our work follows papers of Bressler-Evens \cite{BE1, BE2},
Calm\`es-Petrov-Zainoulline \cite{CPZ}, Harada-Holm-Henriques \cite{HHH}, 
Hornbostel-Kiritchenko \cite{HK}, and Kiritchenko-Krishna \cite{KiKr}, which have laid important foundations.
Combining these tools we study the equivariant cohomology of the flag varieties, partial flag varieties, and Schubert
varieties via the algebraic and combinatorial study of the rings which appear in the Borel model and 
the moment graph model.  In Sections \ref{Borelmodel} and \ref{momgrmodel} we review the setup for these models and
the connection to the (generalized) nil affine Hecke algebra  and the BGG-Demazure operators (see also \cite{SZ} and
\cite{BE1, BE2}). 

One of the main points of our work is to shift the focus from Bott-Samelson classes to Schubert classes. 
In ordinary equivariant cohomology and equivariant K-theory these agree, but in 
generalized  cohomology the Schubert classes and the Bott-Samelson classes usually differ. 
Since the Schubert varieties are not, in general, smooth it is not even clear how the Schubert classes
(the fundamental classes of the Schubert varieties) should be defined.
In Section \ref{Schubertclasses} we give
explicit examples of ``naive pushforwards'' and Bott-Samelson classes and explain why neither of these
can possibly be the Schubert classes in general.  There are several directions to explore in searching for
a good way to define Schubert classes:
\begin{enumerate}
\item[(a)] One can take the lead of Borisov-Libgober \cite{BL} (see also \cite{To}), and define the Schubert class 
$[X_w]$ as a `corrected' version of the Bott-Samelson class $[Z_{\vec w}]$
which, in the end, does not depend on the reduced word $\vec w$ chosen for $w$.  
Borisov-Libgober \cite[Definition 3.1]{BL}
obtain a correction factor for the elliptic genus from the discrepancies of the components of the
exceptional divisor of a resolution of singularities of a variety with at worst log terminal singularities. Recent
papers of Anderson-Stapledon \cite{AS} and Kumar-Schwede \cite{KS} explain that Schubert varieties have
Kawamata log terminal singularities and analyze the exceptional divisor in the 
Bott-Samelson resolution.  In Section \ref{Schubertclasses}
we compute a possible equivariant algebraic cobordism correction factor for the smallest singular 
(complex dimension 3) Schubert variety in all rank 2 cases.
Though the approach of Borisov-Libgober was a motivation for our computations
we have not yet understood how to make our computation of the correction factor for equivariant algebraic
cobordism relate to the correction suggested by Borisov-Libgober for the elliptic genus.
\item[(b)] One can try to define the Schubert classes as classes determined, hopefully uniquely, by positivity
properties under multiplication.  We have not yet managed to make a definition that is satisfying
but our computations of Schubert products do display remarkable positivity features.
\item[(c)] One can try to use the theory of Soergel bimodules (see \cite{Soe})
to pick out particular generators (as $(S,S)$-bimodules)
of the generalized cohomologies of Schubert varieties which serve as Schubert classes.
Though we have not had space to exhibit our computations of the algebraic cobordism case of Soergel bimodules
in this paper, our preliminary computations show that generalizing the Soergel bimodule theory to the ring
$S$ which appears in Theorem \ref{GBmom} is useful for obtaining better understanding of the equivariant
generalized cohomology of Schubert varieties.
\end{enumerate}

In Section \ref{rank2comps} we provide explicit computations of Schubert classes, and 
products with Schubert classes in the rank 2 cases.  
Our computations hold for all rank two cases, but we have only given specific results for
Schubert classes of Schubert varieties in $G/B$ of (complex) dimension $\le 3$.  In partiuclar, this
provides complete results for types $A_2$ and $B_2$ and partial results for $G_2$ and $A_1^{(1)}$.

To some extent this paper is a sequel to \cite{GR}.  That paper considers the case of equivariant
$K$-theory.  In retrospect, \cite{GR} did not capitalize on the full power of
the moment graph model, in particular, that the map $\Phi$ in Theorem \ref{GBmom} is a \emph{ring homomorphism}.  
This key point is the feature which we exploit in this paper to execute computations similar to those in \cite{GR},
but with greater ease and in greater generality.

\medskip\noindent
\textbf{Acknowledgments.}  We thank the Australian Research Council for continuing support of our research
under grants  DP0986774, DP120101942 and DP1095815. Many thanks to Geordie Williamson, Omar Ortiz, and Martina Lanini 
for teaching us the theory of moment graphs and this beautiful way
of working with $T$-equivariant cohomology theories.  We thank Alex Ghitza, Matthew
Ando, Megumi Harada, Dave Anderson and Michel Brion for helpful conversations.  We also thank
Craig Westerland for answering many many questions of all shapes and sizes
all along the way. It is a pleasure to dedicate this paper to C.S. Seshadri who, for so many years, has provided
so much Schubert calculus support and inspiration.

\section{The Schubert calculus framework}\label{flagvars}

\subsection{Flag and Schubert varieties}

The basic data is 
\begin{equation}\label{basicGdata}
\begin{array}{ll}
G &\hbox{a connected complex reductive algebraic group} \\
\cup\vert \\
B &\hbox{a Borel subgroup} \\
\cup\vert \\
T &\hbox{a maximal torus.}
\end{array}
\end{equation}
The \emph{Weyl group}, the \emph{character lattice} and \emph{cocharacter lattice} are, respectively,
\begin{equation}\label{basicWdata}
W_0=N(T)/T, \qquad
\fh_\ZZ^* = \Hom(T, \CC^\times)
\qquad\hbox{and}\qquad
\fh_\ZZ = \Hom(\CC^\times, T),
\end{equation}
where $\Hom(H,K)$ is the abelian group of algebraic group homomorphisms from $H$ to $K$
with product given by pointwise multiplication, $(\phi\psi)(h) = \phi(h)\psi(h)$.  Since the Weyl
group acts on $T$, it also acts on $\fh_\ZZ^*$ and on $\fh_\ZZ$. 

A \emph{standard parabolic subgroup} of $G$ is a subgroup $P_J\supseteq B$ such that
$G/P_J$ is a projective variety. A \emph{parabolic subgroup} of $G$ is a conjugate of a
standard parabolic subgroup.
\begin{equation}\label{flagvardefn}
\hbox{The \emph{flag variety} is $G/B$}
\qquad\hbox{and}\qquad
\hbox{$G/P_J$ are the \emph{partial flag varieties}.}
\end{equation}
These are  studied via the \emph{Bruhat decomposition}
\begin{equation}\label{Bruhatdec}
G = \bigsqcup_{w\in W_0} BwB
\qquad\hbox{and}\qquad
G = \bigsqcup_{u\in W^J} BuP_J
\end{equation}
where
$W_J = \{ v\in W_0\ |\ vT\subseteq P_J\}$ and
\begin{equation}\label{WJ}
W^J = \{ \hbox{coset representatives $u$ of cosets in $W_0/W_J$}\}.
\end{equation}
The \emph{Schubert varieties} are
\begin{equation}\label{Schubertvardefn}
X_w = \overline{BwB}\quad  \hbox{in $G/B$} \qquad\hbox{and}\qquad
X_u^J = \overline{BuP_J}\quad \hbox{in $G/P_J$},
\end{equation}
and the \emph{Bruhat orders}  are the partial orders on $W_0$ and $W_J$ given by
\begin{equation}\label{Bruhat order}
X_w = \overline{BwB} = \bigsqcup_{v\le w} BvB
\qquad\hbox{and}\qquad
X_u^J = \overline{BuP_J} = \bigsqcup_{z\le u} BzP_J.
\end{equation}
The $T$-fixed points
\begin{equation}\label{Tfixpt}
\hbox{in $G/B$ are}\ \{wB\ |\ w\in W_0\} \quad\hbox{and} \qquad
\hbox{in $G/P_J$ are}\ \{uP_J\ |\ u\in W^J\}.
\end{equation}

Let $P_1,\ldots, P_n$ be the minimal parabolic subgroups $P_i\ne B$.  Then
\begin{equation}\label{simplerefl}
W_i = W_{\{i\}} = \{ 1, s_i\}
\qquad\hbox{and}\qquad
\hbox{$s_1,\ldots, s_n$ are the \emph{simple reflections} in $W_0$.}
\end{equation}
With respect to the action of $W_0$ on 
$\fh_\RR^* = \RR\otimes_\ZZ \fh_\ZZ^*$, the $s_i$ are reflections in the 
hyperplanes $(\fh^*)^{s_i}=\{ \mu\in \fh_\RR^*\ |\ s_i\mu = \mu\}$.
An alternative description of the standard parabolic subgroups is to 
let $J\subseteq \{ 1, 2, \ldots, n\}$ and let
\begin{equation}\label{PJdecomp}
W_J = \langle s_j\ |\ j\in J\rangle.
\qquad\hbox{Then}\qquad
P_J = \bigsqcup_{v\in W_J} BvB.
\end{equation}
In particular, $P_i = P_{\{i\}} = B \sqcup Bs_iB$, for $i=1, 2,\ldots, n$.

\begin{thm} (Coxeter)  The group $W_0$ is generated by $s_1,\dots, s_n$ with relations
$$s_i^2 = 1
\qquad\hbox{and}\qquad
\underbrace{s_is_js_i\cdots}_{m_{ij}\ \mathrm{factors}}
= \underbrace{s_js_is_j\cdots}_{m_{ij}\ \mathrm{factors}}
$$
where $\pi/m_{ij} = (\fh^*)^{s_i} \angle (\fh^*)^{s_j}$ is the angle between
$(\fh^*)^{s_i}$ and $(\fh^*)^{s_j}$.
\end{thm}


The definitions in \eqref{flagvardefn}, \eqref{Tfixpt} and \eqref{Schubertvardefn} provide $T$-equivariant maps
\begin{align}\label{flgmorphisms}
\begin{matrix}
p_J\colon &G/B &\longrightarrow &G/P_J \\
&gB &\longmapsto &gP_J
\end{matrix}
\qquad\quad
\begin{matrix}
\iota_w\colon &\pt &\hookrightarrow &G/B \\
&\pt &\longmapsto &wB
\end{matrix}
\qquad\quad
\begin{matrix}
\sigma_w\colon &X_w &\hookrightarrow &G/B \\
&gB &\longmapsto &gB
\end{matrix}
\end{align}
and
\begin{align}\label{ptlflgmorphisms}
\begin{matrix}
\iota_u^J\colon &\pt &\hookrightarrow &G/P_J \\
&\pt &\longmapsto &uP_J
\end{matrix}
\qquad\quad
\begin{matrix}
\sigma_u^J\colon &X_u^J &\hookrightarrow &G/P_J \\
&gP_J &\longmapsto &gP_J
\end{matrix}
\end{align}
for $J\subseteq \{1,2,\ldots, \ell\}$, $w\in W_0$, and $u\in W^J$.

For example, in type $G=GL_3$, with $T$ and $B$ the subgroups given by
$$T =  \left\{\begin{pmatrix} * &0 &0 \\ 0 &* &0 \\ 0 &0 &*\end{pmatrix}\right\}
\qquad\hbox{and}\qquad
B =  \left\{\begin{pmatrix} * &* &* \\ 0 &* &* \\ 0 &0 &*\end{pmatrix}\right\},
$$
then $W_0 = \langle s_1, s_2\ |\ s_i^2=1, s_1s_2s_1=s_2s_1s_2\rangle$, where
$$
s_1 =\begin{pmatrix} 0 &1 &0 \\ 1 &0 &0 \\ 0 &0 &1\end{pmatrix}
\qquad\hbox{and}\qquad
s_2 =\begin{pmatrix} 0 &1 &0 \\ 1 &0 &0 \\ 0 &0 &1\end{pmatrix}.
$$
Then
$$  
  
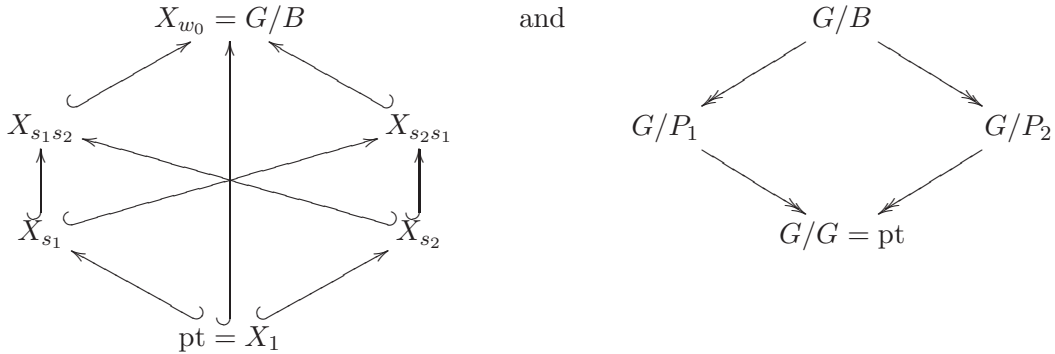
\xymatrix{
  &X_{w_0}=G/B \\
 X_{s_1s_2} \ar@{^{(}->}[ur] & &X_{s_2s_1} \ar@{_{(}->}[ul] \\
 X_{s_1} \ar@{^{(}->}[u] \ar@{^{(}->}[urr] & &X_{s_2} \ar@{^{(}->}[u] \ar@{_{(}->}[ull] \\
 &\pt=X_1 \ar@{^{(}->}[ur] \ar@{_{(}->}[ul] \ar@{^{(}->}[uuu]
 }
 \qquad\hbox{and}\qquad
 \xymatrix{
  &G/B \ar@{->>}[dr] \ar@{->>}[dl] \\
G/P_1\ar@{->>}[dr]  & &G/P_2 \ar@{->>}[dl] \\
 &G/G=\pt
}
$$
where $P_1$ and $P_2$ are the subgroups of $G=GL_3(\CC)$ given by
$$P_1 = \left\{\begin{pmatrix} * &* &* \\ * &* &* \\ 0 &0 &*\end{pmatrix}\right\}
= B \sqcup Bs_1B
\qquad\hbox{and}\qquad
P_2 = \left\{\begin{pmatrix} * &* &* \\ 0 &* &* \\ 0 &* &*\end{pmatrix}\right\}
=B\sqcup Bs_2B.$$

\subsection{Generalized cohomology theories}\label{cohomtheories}

Schubert calculus is the study of the cohomology of flag and Schubert varieties.
Although the home for our computations is the particular ring $S = \LL[[y_\lambda]]$
of \eqref{Rdefn} the motivation comes from the formalism of 
generalized cohomology theories $h$.
Model examples are: ordinary cohomology $H$, $K$-theory $K$,
elliptic cohomology (see \cite{MR, GKV, Gr, An, Lu})
and complex and algebraic cobordism $\Omega$ (see \cite{LM}).
Key to our point of view is that if $f\colon X\to Y$ is a morphism of spaces,
the contravariance of the cohomology theory provides
$$\hbox{a pullback}\quad
f^*\colon h(Y)\to h(X),
\qquad\hbox{and}\quad
\hbox{a pushforward}\quad
f_!\colon h(X)\to h(Y)
$$
exists if the morphism $f$ is nice enough.
Our true interest is in the morphisms in \eqref{flgmorphisms} and \eqref{ptlflgmorphisms}
and \eqref{BSresolution}.
Sometimes we will try to consider, by combinatorial gadgetry, pushforwards across these
morphisms even in cases where we are not sure
that, for any given cohomology theory, the pushforward properly exists.

As in \cite[\S8.2]{CPZ}, the important property for the analysis of Schubert calculus
is that an oriented cohomology theory $h$ comes
with a formal group law $F$ over the coefficient ring $h(pt)$ such that
$$F(c_1^h(\cL_1),c_1^h(\cL_2)) = c_1^h(\cL_1\otimes \cL_2),$$
where $\cL_1$ and $\cL_2$ are line bundles on $X$ and $c_1^h$ denotes the first Chern class
in the cohomology theory $h$ (see \cite[Cor. 4.1.8]{LM}).
The \emph{Lazard ring}
$\LL$ is generated by symbols $a_{ij}$, for $i,j\in \ZZ_{>0}$, which satisfy the relations
given by the equations
\begin{equation}\label{fmgplaw}
F(x,F(y,z))=F(F(x,y),z), \qquad F(x,y)=F(y,x), \qquad F(x,0)=x,
\end{equation}
where
$$F(x,y) = x+y+a_{11}xy+a_{12}xy^2+a_{21}x^2y+\cdots$$
The ring $\LL$ is the universal coefficient ring for a formal group law $F$.  This ring is one of the ingredients
for the construction of the ring $S$ where we do our computations.

A equivariant cohomology theory $h_T$ is a functor from $T$-spaces (some appropriate class of topological
or geometric objects with $T$-action) to some class of algebraic objects (in most of our model examples, $h_T(\pt)$-algebras).
Important features and properties of the theory include:
\begin{enumerate}
\item[(0)] Normalization: specification of $h_T(\pt)$,
\item[(1)] nice behaviour under products, smashes, suspensions: such as axioms for computing $h_{G\times K}(M\times N)$,
\item[(2)] functoriality/pullbacks:
if $f\colon X\to Y$ then we have $f^*\colon h_T(Y)\to h_T(X)$
\item[(3)] Thom isomorphism/orientability/pushforwards: 
For certain classes of maps $f\colon X\to Y$ there exists a pushforward $f_!\colon h_T(X)\to h_T(Y)$,
\item[(4)] Change of groups: For certain classes of groups $G$ and $K$ and group homomorphisms 
$\varphi\colon G\to K$ there exist $\chi_\varphi\colon h_G\to h_K$ and $\chi^\varphi\colon h_K\to h_G$.
\end{enumerate}
The art of choosing appropriate categories of input ``$T$-spaces'', of output ``algebraic objects''
and widening the classes of maps on which pushforwards and/or change of groups homomorphisms
are defined is a beautiful chapter in algebraic topology and geometry. The challenge of
extending a nonequivariant generalized cohomology theory to the equivariant case can be considerable.
For such a genuinely equivariant theory the formal groups above will be replaced by actual groups but we
do not emphasize this point of view here.
For a small selection of references we refer the reader to \cite[p.\ 37-29]{Ad}
for a discussion of the connection to formal group laws and spectra, \cite[Chapt.\ XIII]{Ma} and \cite{Oko}
for a discussion of equivariant orientable theories as Mackey functors and \cite[(1.5)]{GKV}
for discussion of axioms for equivariant elliptic cohomology.
 
In order to specify a home for our computations in Schubert calculus in equivariant cohomology theories we
follow \cite{HHH}.  They 
restrict their class of spaces to \emph{GKM spaces}:
stratified $T$-spaces
$$X = \bigcup_{i\in \ZZ_{>0}} X_i,
\qquad X_1\subseteq X_2\subseteq X_3\subseteq \cdots,
$$
where the successive quotients $X_i/X_{i-1}$ are homeomorphic to the Thom
spaces $Th(V_i)$ of some $h$-orientable $T$-vector bundles $V_i\to F_i$ (see \cite[(2.1)]{HHH}).
As pointed out in \cite[Remark 3.3]{HHH}, for the case of flag and Schubert varieties that
are the focus of this paper, the $F_i$ are points and the $V_i$ are one dimensional
representations of $T$.  In particular, the assumptions of \cite[\S 3]{HHH} hold for these cases.

\subsection{The Borel model for $h_T(G/B)$}\label{Borelmodel}

The general combinatorial Schubert calculus uses $\fh_\ZZ^*$ and $\fh_\ZZ$ to
build a $C$-algebra $R$ with an action of $W_0$ on $R$ by $C$-algebra automorphisms
(in favorite examples $C$ may be $\ZZ$, or the ring $\widetilde{Th}_0$ of holomorphic
functions on the upper half plane, or the Lazard ring $\LL$, see the examples below).
If 
$$R^{W_0} = \{ f\in R\ |\ \hbox{$wf=f$ for $w\in W_0$}\}
\qquad\hbox{is the \emph{invariant ring},}$$
then, conceptually,
\begin{equation}
R = h_T(\pt)
\qquad\hbox{and}\qquad
R^{W_0} = h_G(\pt),
\end{equation}
for the equivariant cohomology theory $h_T$ under analysis.  
By definition, the \emph{coinvariant ring} is
\begin{equation}\label{coinvdefn}
R\otimes_{R^{W_0}}R
=
\frac{R\otimes_C R}{\langle f\otimes 1 - 1\otimes f \ |\ f\in R^{W_0} \rangle},
\end{equation}
where the terminology is chosen to be representative of the classical terminology
in the study of the cohomology of $G/B$, not to reflect a notion of coinvariants with respect
to a group action.
Then (see \cite[Proposition 26.1]{Bo}, \cite[Proposition 1.6]{KL}, \cite[Theorem 4.7]{KiKr}) the ring
\begin{equation}\label{combmodel}
R\otimes_{R^{W_0}}R
\qquad\hbox{is a good combinatorial model for}\quad h_T(G/B),
\end{equation}
where the product on $R\otimes_{R^{W_0}}R$ is given by 
$(f_1\otimes g_1)(f_2\otimes g_2) = f_1f_2\otimes g_1g_2$.

There are four favorite examples:

\medskip\noindent
\textbf{Cohomology: $h_T = H_T$.}  Here
$$H_T(\pt) = S(\fh_\ZZ^*) = \CC[x_1,\ldots, x_n]
\qquad\hbox{and}\qquad
H_G(\pt) = H_T(\pt)^{W_0} = \CC[x_1,\ldots, x_n]^{W_0},
$$
where $x_i = x_{\omega_i}$, where
$\omega_1,\ldots, \omega_n$ is a $\ZZ$-basis of $\fh_\ZZ^*$.
Alternatively, $H_T(\pt)$ is the ring
$$\CC[x_\lambda\ |\ \lambda\in \fh_\ZZ^*]
\quad\hbox{with}\quad
x_{\lambda+\mu} = x_\lambda + x_\mu,
$$
for $\lambda,\mu\in \fh_\ZZ^*$ and with $wx_\lambda = x_{w\lambda}$ for $w\in W_0$ and $\lambda\in \fh_\ZZ^*$.
Then
\begin{align*}
H_T(G/B) 
= H_T(\pt)\otimes_{H_G(\pt)} H_T(\pt) 
= \frac{\CC[y_1,\ldots, y_n, x_1,\ldots, x_n]}
{\langle f(x_1,\ldots, x_n)-f(y_1,\ldots, y_n)\ |\ 
f\in \CC[x_1,\ldots, x_n]^{W_0}\rangle}.
\end{align*}

\medskip\noindent
\textbf{K-theory: $h_T=K_T$.}  Here
$$K_T(\pt) = \CC[\fh_\ZZ^*] = \CC[X_1^{\pm1},\ldots, X_n^{\pm1}]
\qquad\hbox{and}\qquad
K_G(\pt) = K_T(\pt)^{W_0} = \CC[X_1^{\pm1},\ldots, X_n^{\pm1}]^{W_0},
$$
where $X_i = e^{\omega_i}$, where
$\omega_1,\ldots, \omega_n$ is a $\ZZ$-basis of $\fh_\ZZ^*$.
Alternatively, $K_T(\pt)$ is the ring
$$\CC[e^\lambda\ |\ \lambda\in \fh_\ZZ^*]
\quad\hbox{with}\quad
e^{\lambda+\mu} = e^\lambda e^\mu,
$$
for $\lambda,\mu\in \fh_\ZZ^*$ and with $we^\lambda = e^{w\lambda}$ for $w\in W_0$ and $\lambda\in \fh_\ZZ^*$.
Then
\begin{align*}
K_T(G/B) = K_T(\pt)\otimes_{K_G(\pt)} K_T(\pt) 
= \frac{\CC[Y_1^{\pm1},\ldots, Y_n^{\pm1}, X_1^{\pm1},\ldots, X_n^{\pm1}]}
{\langle f(X_1,\ldots, X_n)-f(Y_1,\ldots, Y_n)\ |\ 
f\in \CC[X_1^{\pm1},\ldots, X_n^{\pm1}]^{W_0}\rangle}.
\end{align*}

\medskip\noindent
\textbf{Elliptic cohomology: $h_T = Ell_T$.}  Here $Ell_T(\pt)$ is the structure sheaf of the abelian
variety $A_\tau = \fh_\CC^*/(\fh_\ZZ^*+\tau\fh_\ZZ^*)$.  The
homogeneous coordinate ring
$$\hbox{for $A_\tau$ is}\qquad
\widetilde{Th} = \bigoplus_{m\in \ZZ_{\ge 0}} \widetilde{Th}_m,
\qquad\hbox{and}\qquad
\hbox{for $A_\tau/W_0$ is}\qquad \widetilde{Th}^{W_0}.
$$
Then the graded $\widetilde{Th}$-module corresponding to
$$\hbox{the sheaf
$Ell_T(G/B)$ on $A_\tau$ is}
\qquad
\widetilde{Th}\otimes_{\widetilde{Th}^{W_0}} \widetilde{Th}.
$$

\medskip\noindent
\textbf{Complex or algebraic cobordism: $h_T = \Omega_T$.}  Algebraic cobordism is treated in the book of
Levine-Morel \cite{LM} and $T$-equivariant algebraic cobordism $\Omega_T$ is treated in 
\cite{Kr} and \cite{KiKr}.  The following summary of our setting
is made precise by Theorem \ref{GBmom} below.

The \emph{Lazard ring} $\LL$ is the coefficient ring for the universal formal group
law $F$ so that $\LL$ is given by generators $a_{ij}$ with relations given by setting
$$F(x,y) = x+y+\sum_{i,j\in \ZZ_{>0}} a_{ij}x^iy^j
\qquad\hbox{in $\LL[[x,y]]$,}
$$
and requiring
$$F(x,0)=F(0,x)=x, \quad
F(x,y)=F(y,x), \quad
F(x, F(y,z))=F(F(x,y),z).$$
Then
$$\Omega_T(\pt) = 
\LL[[x_\lambda\ |\ \lambda\in \fh_\ZZ^*]]
\quad\hbox{with}\quad
x_{\lambda+\mu} = x_\lambda +_F x_\mu = F(x_\lambda,x_\mu),
$$
for $\lambda,\mu\in \fh_\ZZ^*$. Then
$$\Omega_G(\pt) = \Omega_T(\pt)^{W_0} = \LL[[x_\lambda\ |\ \lambda\in \fh_\ZZ^*]]^{W_0},
\qquad\hbox{where}\quad
wx_\lambda = x_{w\lambda},
$$
for $w\in W_0$ and $\lambda\in \fh_\ZZ^*$,
and
\begin{align*}
\Omega_T(G/B) 
= \Omega_T(\pt)\otimes_{\Omega_G(\pt)} \Omega_T(\pt) 
= \frac{\LL[[y_\lambda, x_\mu\ |\ \lambda\in \fh_\ZZ^*]]}
{\langle f(x)-f(y)\ |\ 
f\in \LL[[x_\lambda\ |\ \lambda\in \fh_\ZZ^*]]^{W_0}\rangle}.
\end{align*}

\medskip\noindent
Sample references for such identities are \cite{KK1} for the case of $H_T(G/B)$, \cite{KK2,KL,CG} for $K_T(G/B)$,
\cite{KP,Gr, GKV, An, Ga} for $Ell_T(G/B)$ and \cite{HHH, CPZ, HK, KiKr} for $\Omega_T(G/B)$.

The cobordism case specializes to the cases of cohomology $H_T$ and $K$-theory $K_T$ by setting
$$
F(x,y) = \begin{cases}
x+y, &\hbox{in $H_T$,} \\
x+y-xy, &\hbox{in $K_T$,}
\end{cases}
\qquad\hbox{and}\qquad
x_\lambda = \begin{cases}
x_\lambda, &\hbox{in $H_T$,} \\
1-e^\lambda, &\hbox{in $K_T$.}
\end{cases}
$$

\section{The moment graph model}\label{momgrmodel}

\subsection{$T$-fixed points and the map $\Phi$}

Following Goresky-Kottwitz-MacPherson \cite[Theorem 1.2.2]{GKM} a powerful
way to think about this theory is via the \emph{moment graph model}.  This means that
for a $T$-variety $X$ where the imbeddings of the $T$-fixed points of $X$ into $X$ are
\begin{equation}
\begin{matrix}
\iota_w\colon &\pt &\to &X \\
&* &\mapsto &w
\end{matrix}
\qquad\hbox{consider}\qquad
\hbox{$\iota^*=\bigoplus_{w\in W} \iota_w^*$}\colon \Omega_T^*(X) \mapright{} \bigoplus_{w\in W} \Omega_T(\pt),
\end{equation}
where the sums are over an index set $W$ for the $T$-fixed points in $X$.  When $X$ is a 
``GKM-space" (see \cite[Theorem 14]{GKM} for several equivalent characterization of a GKM space for
equivariant ordinary cohomology and \cite[HHH] for equivariant generalized cohomology theories)
the ring homomorphism $\iota^*$ is injective with image
$$
\im\,\iota^* = \left\{ (g_w)_{w\in W_0} \in \bigoplus_{w\in W_0} \Omega_T(\pt),\ \bigg|\ 
\begin{matrix}
\hbox{$g_w-g_{w'} \in y_\alpha \Omega_T(\pt)$ if there is a} \\
\hbox{1-dimensional $T$-orbit containing $w$ and $w'$}
\end{matrix} \right\},
$$
where $y_\alpha$ is the $T$-equivariant Chern class of the tangent along the 1-dimensional
orbit connecting $w$ and $w'$.

Computations are facilitated by encoding the information of $\im\, \iota^*$ with a 
\emph{moment graph}, which has vertices corresponding to the $T$-fixed points of $X$
and labeled edges $w\mapright{\alpha}w'$ corresponding to 1-dimensional $T$-orbits in $X$.  
For example,
for $G/B$ for type $GL_3$ the graph is
\begin{equation}\label{momgrA2}
\xymatrix{
  &1 \ar[dl]_{y_{-\alpha_1}} \ar[dr]^{y_{-\alpha_2}} \ar[ddd]^<<<<<<<<<{y_{-(\alpha_1+\alpha_2)}} \\
s_1 \ar[d]_{y_{-\alpha_2}} \ar[drr]_<<<<<<<<<<{y_{-(\alpha_1+\alpha_2)}} & &s_2\ar[d]^{y_{-\alpha_1}}
\ar[dll]^<<<<<<<<<{y_{-(\alpha_1+\alpha_2)} } \\
s_1s_2 \ar[dr]_{y_{-\alpha_1}} & &s_2s_1 \ar[dl]^{y_{-\alpha_2}} \\
 &s_1s_2s_1=s_2s_1s_2
 }
\end{equation}
A \emph{moment graph section} is a tuple $(g_w)_{w_\in W}$ of elements of $\Omega_T(\pt)$
which is an element of $\im\,\iota^*$.

A morphism of GKM-spaces is a morphism of $T$-spaces
$$f\colon X\to Y
\qquad\hbox{which provides, by restriction,}\qquad
f\colon W\to V$$
from the set $W$ of $T$-fixed points of $X$ to the 
set $V$ of $T$-fixed points of $Y$.
Viewing elements of $H_T(X)$ and $H_T(Y)$ as moment graph sections
the maps
$$f^*\colon H_T(Y)\to H_T(X)
\qquad\hbox{and}\qquad
f_!\colon H_T(X)\to H_T(Y)
$$
are given by
\begin{equation}\label{pullandpush}
(f^*(c))_w = c_{f(w)}, \qquad\hbox{and}\qquad
(f_!(\gamma))_v = \sum_{w\in f^{-1}(v)} \gamma_w\frac{1}{e(f)_{wv}},
\end{equation}
where the \emph{Euler class} of $f$ from $v$ to $w$ is 
$$e(f)_{wv} = \left(\prod_{{\mathrm{edges\ of}\ W \atop \mathrm{adjacent\ to}\ w} } y_\beta\right)
\left(\prod_{{\mathrm{edges\ of}\ V \atop \mathrm{adjacent\ to}\ v} } y_\beta\right)^{-1}.$$
The second formula in \eqref{pullandpush} is a form of the 
familiar formula for push forwards by ``localization at the $T$-fixed points'' as
found, for example, in \cite[(3.8)]{AB}.  The Euler class of $f$ from $v$ to $w$ is the contribution
measured by the difference between 
the tangent space at the $T$-fixed point $w$ in $X$ to the tangent space to the $T$-fixed point $v=f(w)$ in $Y$.

The Borel model and the moment graph model for $G/B$ for
equivariant algebraic cobordism $\Omega_T(G/B)$ are
summarized in the following Theorem, which is a combination of  
\cite[Theorem 4.7]{KiKr} and \cite[Theorem 3.1]{HHH}.
The ring $S$ which takes the role of $\Omega_T(\pt)$ is as in \cite[\S2.4]{CPZ}.
For comparison to the $K$-theory case see \cite[Theorem 3.13]{KK2} and \cite[Theorem 3.1]{LSS}.

\begin{thm}\label{GBmom} \emph{(\cite[Theorem 3.1]{HHH}, \cite[Theorem 4.7]{KiKr}
and \cite[\S2.4]{CPZ} combined)} 
Let $G\supseteq B\supseteq T$ be a reductive group datum as in \eqref{basicGdata} and 
let $W_0$ and $\fh_\ZZ^*$ be the Weyl group and the weight lattice $\fh_\ZZ^*$ as in
\eqref{basicWdata}.
Let $\LL$ be the Lazard ring generated by $a_{ij}$ as in \eqref{fmgplaw}
and let $S$ be the $\LL$-algebra
\begin{equation}\label{Rdefn}
S = \LL[[y_\lambda\ |\ \lambda\in \fh_\ZZ^*]],
\quad\hbox{with}\quad
y_{\lambda+\mu} =  y_\lambda+y_\mu+a_{11}y_\lambda y_\mu+a_{12}y_\lambda y_\mu^2
+a_{21}y_\lambda^2 y_\mu + \cdots.
\end{equation}
The Weyl group
$$W_0\quad\hbox{acts $\LL$-linearly on $S$ by } \qquad
wy_\lambda = y_{w\lambda},
$$
for $w\in W_0$, $\lambda\in \fh_\ZZ^*$.
Define a product on $\bigoplus_{w\in W_0} S$ pointwise,
\begin{equation}\label{momgrproduct}
(f_w)_{w\in W_0} \cdot (g_w)_{w_\in W_0} = (f_wg_w)_{w\in W_0},
\end{equation}
and let $S\otimes_{S^{W_0}} S$ be the coinvariant ring as defined in \eqref{coinvdefn}.
The $S$-algebra homomorphism
\begin{equation}\label{momgrmap}
\xymatrix{
\Phi\colon S\otimes_{S^{W_0}} S \ar[r]^{\sim} 
&\Omega_T(G/B) \ar[r]^{\sim} &\im\,\Phi \ar@{^{(}->}[r] 
&\bigoplus_{w\in W_0} S \\
\quad\quad f\otimes g\quad \ar@{|->}[rrr] 
&&&\quad\big(f\cdot (w^{-1}g)\big)_{w\in W_0} 
} 
\end{equation}
is well defined and injective with 
$$
\im\,\Phi = \left\{ (g_w)_{w\in W_0} \in \bigoplus_{w\in W_0} S\ \bigg|\ 
\hbox{$g_w-g_{ws_\alpha} \in y_{-\alpha} S$ for $\alpha\in R^+$ and $w\in W_0$} \right\},
$$
where $R^+$ is the set of positive roots corresponding to $B$ and
$s_\alpha\in W_0$ denotes the reflection corresponding to $\alpha$.
\end{thm}

To provide a feel for the ring $S$ of \eqref{Rdefn}, let us provide some formulas which will be useful for
computations later. To recapitulate and summarize previous definitions,
\begin{equation}\label{Sdefn}
S = \LL[[y_\lambda\ |\ \lambda\in \fh_\ZZ^*]]
\qquad\hbox{with}\qquad y_{\lambda+\mu} = y_\lambda+y_\mu-p(y_\lambda,y_\mu)y_\lambda y_\mu,
\end{equation}
where $p(y_\lambda,y_\mu)\in \LL[[y_\lambda, y_\mu]]$ is a power series
\begin{equation}\label{pdefn}
p(y_\lambda, y_\mu) = -a_{11}-a_{12}y_\mu-a_{21}y_\lambda-a_{31}y_\lambda^2-a_{22}y_\lambda y_\mu
-a_{13}y_\mu y_\lambda - \cdots,
\end{equation}
with $a_{ij}\in \LL$ satisfying relations such that
\begin{equation}
y_{-\lambda+\lambda}=y_0=0,\qquad
y_{\lambda+\mu}=y_{\mu+\lambda}, \qquad y_{(\lambda+\mu)+\nu}=y_{\lambda+(\mu+\nu)}.
\end{equation}
Then
\begin{equation}\label{yalphynegalapha}
y_\alpha = \frac{-y_{-\alpha}}{1-p(y_\alpha,y_{-\alpha})y_{-\alpha}},
\qquad 
\frac{1}{y_{-\alpha}}+\frac{1}{y_\alpha}
= p(y_\alpha, y_{-\alpha}),
\end{equation}
and the formula
\begin{equation}\label{yellalphaoveryalpha}
\frac{y_{-\ell\alpha}}{y_{-\alpha}}
= \ell - \sum_{j=1}^{\ell-1} p(y_{-\alpha},y_{-j\alpha})y_{-j\alpha}
= 1 + \sum_{j=1}^{\ell-1} (1-p(y_{-\alpha}, y_{-j\alpha})y_{-j\alpha}),
\qquad\hbox{for $\ell\in \ZZ_{>0}$,}
\end{equation}
is proved by induction on $\ell$.  Using \eqref{yellalphaoveryalpha} and the formula
$s_i\lambda = \lambda - \langle \lambda, \alpha_i^\vee\rangle\alpha_i$ for the action of a 
simple reflection on $\fh^*$ produces
\begin{equation}\label{divdiffexpansion}
\frac{y_{s_i\lambda}-y_{\lambda}}{y_{-\alpha_i}}
=(1-p(y_\lambda,y_{-\langle \lambda, \alpha_i^\vee\rangle\alpha_i})y_\lambda)
\left(1+\sum_{j=1}^{\langle \lambda, \alpha^\vee\rangle - 1}
(1-p(y_{-\alpha_i},y_{-j\alpha_i})y_{-j\alpha_i})\right),
\end{equation}
for $\langle \lambda, \alpha_i^\vee\rangle\in \ZZ_{\ge 0}$.
Formula \eqref{divdiffexpansion}
generalizes one of the favorite formulas for the action of a Demazure operator (see \cite[Lemma 8.2.8]{Ku2}).
This cobordism case specializes to $H_T$ and $K_T$ by setting
\begin{equation}\label{spectoHTandKT}
p(y_\lambda, y_\mu) = \begin{cases}
0, &\hbox{in $H_T$,} \\
1, &\hbox{in $K_T$,}
\end{cases}
\qquad\hbox{and}\qquad
y_\lambda = \begin{cases}
y_\lambda, &\hbox{in $H_T$,} \\
1-e^{\lambda}, &\hbox{in $K_T$.}
\end{cases}
\end{equation}

\subsection{The nil affine Hecke algebra}

Let $S$ be as in \eqref{Rdefn} and \eqref{Sdefn}.  The point of view of \cite{GR} is that the homomorphism
$\Phi$ of \eqref{momgrmap} arises naturally from the nil affine Hecke algebra.

The \emph{nil affine Hecke algebra} $H$ is
\begin{align*}
H &= (S\otimes_\LL S)\ltimes \LL[W_0] \\ 
&= S\hbox{-span}\{ gt_w\ |\ g\in S, w\in W_0\}
= \LL\hbox{-span}\{ (f\otimes g)t_w\ |\ f,g\in S, w\in W_0\} \nonumber
\end{align*}
with
\begin{equation}
\label{naffHeckedefn}
t_ut_v = t_{uv} \quad\hbox{and}\quad t_w(f\otimes g)  = (f\otimes (wg))t_w,
\end{equation}
for $u,v, w\in W_0$ and $f,g\in S$.  
The nil affine Hecke algebra $H$ acts on $S\otimes_\LL S$ and on $S\otimes_{S^{W_0}} S$ by
\begin{equation}\label{nilhactionpoly}
t_w(f\otimes g)= f\otimes wg
\qquad\hbox{and}\qquad (h\otimes p)(f\otimes g) = hf\otimes pg,
\end{equation}
for $h,p,f,g\in S$ and $w\in W_0$.  These actions arise from the realization of $S\otimes_{S^{W_0}} S$
as an induced up $H$-module in \eqref{hTHmod} below.

Let $b_1$ be a symbol and let $Sb_1$ be the $S\otimes_\LL S$ module (a rank
1 free $S$-module with basis $\{ b_1 \}$) 
corresponding to the ring homomorphism
$$\begin{matrix}\label{aughom}
\varepsilon\colon &S\otimes_\LL S &\longrightarrow &S \\
&f\otimes g &\longmapsto &fg
\end{matrix}
\qquad\hbox{so that the $S\otimes_\LL S$ action on $Sb_1$ is given by
$(f\otimes g)b_1 = fgb_1$,}
$$
for $f,g\in S$.  The induced module
$$Hb_1 = \Ind_{S\otimes_\LL S}^{H}(Sb_1)
\qquad\hbox{has $S$-basis}\qquad \{ b_w |\ w\in W_0\},
\quad\hbox{where $b_w = t_wb_1$}.
$$
Let $\mathbf{1}_0 = \sum_{w\in W_0} t_w$.
With the definition of the $H$ action on $S\otimes_\LL S$ as in \eqref{nilhactionpoly},
the sequence of maps (see \cite[Theorem 2.12]{GR})
\begin{equation}\label{hTHmod}
\begin{matrix}
S\otimes_\LL S &\longrightarrow &H\mathbf{1}_0 &\hookrightarrow &H
&\longrightarrow &Hb_1 &\cong \bigoplus_{w\in W_0} S \\
(f\otimes g) &\longmapsto &(f\otimes g)\mathbf{1}_0 \\
&&& &h &\longmapsto &hb_1
\end{matrix}
\end{equation}
is a homomorphism of $H$-modules (with kernel generated by
$\{ f\otimes 1 - 1\otimes f\ |\ f\in S^{W_0}\}$). The maps in \eqref{hTHmod} allow for the 
expansion of any element of $S\otimes_\LL S$ in terms of the basis $\{b_w\ |\ w\in W_0\}$
of $Hb_1$, giving
\begin{align*}
(f\otimes g)\mathbf{1}_0b_1 
&=(f\otimes g)\big(\sum_{w\in W_0} t_w\big)b_1
=\sum_{w\in W_0} t_w (f\otimes(w^{-1}g)) b_1 \\
&=\sum_{w\in W_0} t_w(f\cdot(w^{-1}g)) b_1
=\sum_{w\in W_0} (f\cdot(w^{-1}g)) b_w.
\end{align*}
This formula illustrates that
computing $\Phi(f\otimes g)$ in \eqref{momgrmap} is equivalent to expanding
$(f\otimes g)b_1$ in terms of the $b_w$.
Because of this we use \eqref{momgrmap} and \eqref{hTHmod} to
$$\hbox{identify\quad $\Omega_T(G/B) = Hb_1 = \hbox{$S$-span}\{ b_w\ |\ w\in W_0\} \cong \bigoplus_{w\in W_0} S$}$$
and write elements 
\begin{equation}\label{momgrnotation}
f\in \Omega_T(G/B)
\qquad\hbox{as}\qquad f= \sum_{w\in W_0} f_wb_w.
\end{equation}
The product in $\Omega_T(G/B)$ is then given by \eqref{momgrproduct}.
To more easily keep track of the left and right factors in 
$S\otimes_\LL S$ use the notation
\begin{equation}
x_\mu = 1\otimes y_\mu
\qquad\hbox{and}\qquad
y_\mu = y_\mu\otimes 1.
\end{equation}
Then the formulas
\begin{equation}\label{nilhmomgr1}
x_\lambda\cdot 1 = x_\lambda\sum_{w\in W_0} t_wb_1 = \sum_{w\in W_0} t_w x_{w^{-1}\lambda} b_1
=\sum_{w\in W_0} y_{w^{-1}\lambda}b_w, \quad\hbox{and}
\end{equation}
\begin{equation}\label{nilhmomgr2}
t_v \sum_{w\in W_0} f_w b_w
= \sum_{w\in W_0} f_w t_vb_w
= \sum_{w\in W_0} f_w b_{vw}
= \sum_{z\in W_0} f_{v^{-1}z} b_z,
\end{equation}
provide the formulas for action of the nil affine Hecke algebra in terms of moment graph sections (see \eqref{nilhactionpoly}).
We often view the values $f_w$ as labels on the vertices of the moment graph so that, for exmaple, in type $GL_3$
where the moment graph is as in \eqref{momgrA2},
\eqref{nilhmomgr1} can be written
$$x_\lambda ~~=~~
\begin{matrix}
&y_\lambda \\
y_{s_1\lambda} &&y_{s_2\lambda} \\
y_{s_2s_1\lambda} &&y_{s_1s_2\lambda} \\
&y_{s_1s_2s_1\lambda}
\end{matrix}
$$

\section{Partial flag varieties and Bott-Samelson classes $[Z_{\vec w}]$}\label{BSclasses}

In this section we review the formulas for the Bott-Samelson classes as established in, for example,
\cite{HK, CPZ, BE1, BE2}.  Though some of these references are not considering the equivariant
case, the same machinery applies to define these classes in $\Omega_T(G/B)$.
In particular, this is the place in the theory where the BGG/Demazure operators are derived from
the geometry.  These operators play a fundamental role in the combinatorial study of $\Omega_T(G/B)$.

\subsection{Pushforwards to partial flag varieties: BGG/Demazure operators}

Using the notation for parabolic subgroups and partial flag varieties as
in \eqref{flagvardefn}, if $J\subseteq \{ 1,2,\ldots, n\}$ and 
$$\begin{matrix}
\pi_J\colon &G/B &\to &G/P_J \\
&gB &\mapsto &gP_J
\end{matrix}
\qquad\hbox{then}\qquad
\pi_J(wB) = uP_J, \quad\hbox{where $wW_J=uP_J$}.$$
Then, in the setting of Theorem \ref{GBmom},
$$S\otimes_{S^W_0}S^{W_J}\cong \Omega_T(G/P_J),$$
and $\pi_J^*\colon \Omega_T(G/P_J)\to \Omega_T(G/B)$
and $(\pi_J)_!\colon \Omega_T(G/B)\to \Omega_T(G/P_J)$
 correspond to
\begin{equation}\label{GPpushpul[}
\pi_J^*\colon S\otimes_{S^{W_0}} S^{W_J} \hookrightarrow S\otimes_{S^{W_0}} S
\qquad\hbox{and}\qquad
(\pi_J)_!\colon S\otimes_{S^{W_0}} S \longrightarrow S\otimes_{S^{W_0}} S^{W_J}
\end{equation}
where $(\pi_J)_!$ is given by the operator in the nil affine Hecke algebra given by
$$(\pi_J)_! = \left(\sum_{v\in W_J} t_v\right) \frac{1}{x_J},
\qquad\hbox{where}\quad x_J = \prod_{\alpha\in R^+_J} x_{-\alpha}.
$$
with $R_J^+$ the set of positive roots for $P_J\supseteq B\supseteq T$.
A special case is when $J= \{ i\}$, for which
\begin{equation}\label{pushpull}
W_J = \{1, s_i\}
\quad\hbox{and}\quad
\pi_i^*(\pi_i)_! = A_i = (1+t_{s_i})\frac{1}{x_{-\alpha_i}},
\end{equation}
is the \emph{BGG-Demazure operator} (see \cite[Cor.-Def.\ 1.9]{BE1}).
The calculus of the operators $A_i$ is controlled via the identities in Section \ref{BGGopcalc}.

\subsection{Bott-Samelson classes}

For a sequence $\vec w = (i_1,\ldots, i_\ell)$ with $1\le i_1,\ldots, i_\ell\le n$ define the
\emph{Bott-Samelson class}
\begin{equation}\label{BSclasses}
[Z_{\vec w}] = [Z_{i_1i_2\cdots i_\ell}] = A_{i_1}A_{i_2}\cdots A_{i_\ell}[Z_\pt],
\end{equation}
where, in the notation of \eqref{momgrnotation},
\begin{equation}\label{Zptdefn}
[Z_\pt]_v = \begin{cases}
\prod_{\alpha\in R^+} y_{-\alpha}, &\hbox{if $v=1$}, \\
0, &\hbox{if $v\ne 1$.}
\end{cases}
\end{equation}

\begin{thm} \emph{(\cite[Prop. 1]{BE2}, \cite[Prop. 3.1]{HK}, \cite[Lemma 3.15]{KK2}, see also  \cite[Proposition 4.1]{HHH}) }
The generalized cohomology 
$$\hbox{$h_T(G/B)$ has $h_T(\pt)$-basis}
\qquad
\{ [Z_{\vec w}]=[\gamma_{\vec w}\colon \Gamma_{\vec w}\to G/B]\ |\ w\in W_0\},
$$
where, for each $w\in W_0$, $\vec w=s_{i_1}\cdots s_{i_\ell}$ is  a fixed reduced word for $w$.
\end{thm}

\noindent
Let us explain where this comes from.
Let $X$ be a $T$-variety.  Following \cite[Example 1.9.1]{Fu}, or \cite[\S 5.5]{CG},
a \emph{cellular decomposition} of $X$ is a filtration
$$\emptyset=X_{-1}\subseteq X_0 \subseteq X_1 \subseteq \cdots \subseteq X_d=X$$
by closed subvarieties such that $X_i=X_{i-1}$ are isomorphic to a disjoint union of affine spaces
$\AA^{\ell_i}$ for $i=1,2,\ldots, d$.  The ``cells'' of $X$ are the $X_i-X_{i-1}$.

\begin{thm} \emph{(see \cite[Prop.\ 7]{G}; \cite[Example 1.9.1]{Fu} who refers to \cite{Ch};
\cite[Lemma 5.5.1]{CG}; \cite[Proposition 1]{BE2}; \cite[Theorem 2.5]{HK})}
Let $X$ be a $T$-variety with a cellular decomposition.
Then $h_T(X)$ has an $h_T(\pt)$-basis
given by resolutions of cell closures (choose one
resolution for each cell).
\end{thm}

\noindent
For $X=G/B$, the Bruhat decomposition
$$G = \bigsqcup_{w\in W_0} BwB
\qquad\hbox{provides the desired cell decomposition}
$$
and the \emph{Schubert varieties} $X_w=\overline{BwB}$ are the closures of the 
Schubert cells.  Let $P_1,\ldots, P_n$ be the minimal parabolics of $G$ (with $P_i\supseteq B$ 
and $P_i\ne B$) and let $s_1,\ldots, s_n$ be the corresponding \emph{simple reflections}
in $W_0$.  The group $W_0$ is generated by $s_1,\ldots, s_n$.
Let $\vec w=s_{i_1}\cdots s_{i_\ell}$ be a reduced word for $w$.  Then
the Bott-Samelson variety $\Gamma_{i_1,\ldots, i_\ell} = P_{i_1}\times_B P_{i_2}\times_B \cdots \times_B
P_{i_\ell}/B$ provides a resolution of $X_w$, 
\begin{equation}\label{BSresolution}
\begin{matrix}
\gamma_{i_1,\ldots, i_\ell}\colon &P_{i_1}\times_B P_{i_2}\times_B \cdots \times_B
P_{i_\ell}\times_B \pt &\longrightarrow &X_w &\hookrightarrow G/B \\
&[g_1,\ldots, g_\ell] &\longmapsto &g_1\cdots g_\ell B
\end{matrix}
\end{equation}
Then following, for example, the proof of \cite[Prop.\ 2]{BE2}, since the diagram
\begin{equation}\label{BEdiagram}
\xymatrix{
P_{i_1}\times_B \cdots \times_B P_{i_\ell}\times_B P_{i_{\ell+1}}\times_B\pt
\ar^>>>>>>>>>>>>>>>{\gamma_{i_1\ldots i_{\ell+1}}}[rr] \ar[d]_{\tau}&&G/B \ar^{\pi_{i_{\ell+1}}}[d] \\
P_{i_1}\times_B \cdots \times_B P_{i_\ell}\times_B\pt 
 \ar[r]_>>>>>>>{\gamma_{i_1\ldots i_\ell}} &G/B \ar_{\pi_{i_{\ell+1}}}[r]&G/P_{i_{\ell+1}}
}
\end{equation}
\begin{enumerate}
\item[(a)] commutes, and
\item[(b)] has both vertical maps fibrations with fibre $P_{i_{\ell+1}}/B$,
\end{enumerate}
it is a pullback square.  Thus
\begin{align}
(\gamma_{i_1\ldots i_{\ell+1}})_!(\iota^*(1)) 
 &= \pi_{i_{\ell+1}}^* (\pi_{i_{\ell+1}}\circ\gamma_{i_1\ldots i_\ell})_!(1) \nonumber \\
&= \pi_{i_{\ell+1}}^* (\pi_{i_{\ell+1}})_! (\gamma_{i_1\ldots i_\ell})_!(1) 
= A_{i_{\ell+1}} (\gamma_{i_1\ldots i_\ell})_!(1).
\label{BGGsquare}
\end{align}
The following result then follows by induction.

\begin{thm} \label{BSclassthm} \emph{(\cite[Theorem 3.2]{HK}, \cite[Proposition 2]{BE2})} 
If $I= (i_1,\ldots, i_\ell)$ is a sequence in $\{1,\ldots, n\}$ and
$\gamma_{i_1\ldots i_\ell}$ is as in \eqref{BSresolution} then 
$$[Z_{i_1\cdots i_\ell}]=[(\gamma_{i_1\ldots i_\ell})_!(1)] = A_{i_1}\cdots A_{i_\ell}[Z_\pt],
\qquad\hbox{where $[Z_{\pt}]$ is the class of a point.}
$$
\end{thm}

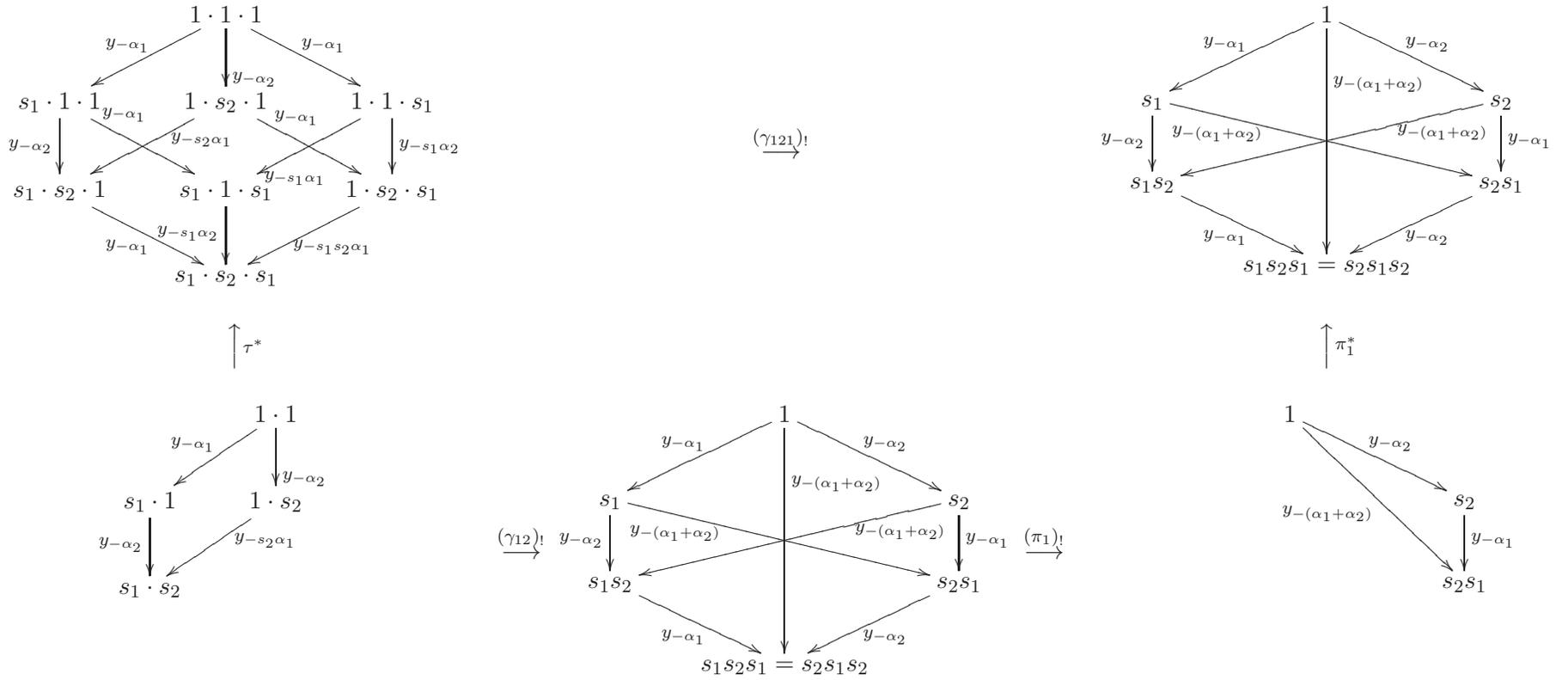
\begin{sidewaysfigure}
$$
\begin{matrix}
  \xymatrix{
  &1\cdot 1\cdot 1 \ar[dl]_{y_{-\alpha_1}} \ar[dr]^{y_{-\alpha_1}} \ar[d]^>>{y_{-\alpha_2}} \\
s_1\cdot 1\cdot 1 \ar[d]_{y_{-\alpha_2}} \ar[dr]^<<<<{y_{-\alpha_1}} 
&1\cdot s_2\cdot 1 \ar[dl]^<<<{y_{-s_2\alpha_1}} \ar[dr]^<<<<<{y_{-\alpha_1}} 
&1\cdot 1\cdot s_1 \ar[d]^{y_{-s_1\alpha_2}} \ar[dl]^>>>>{y_{-s_1\alpha_1}} \\
s_1\cdot s_2\cdot 1 \ar[dr]_{y_{-\alpha_1}} 
&s_1\cdot 1\cdot s_1 \ar[d]_{y_{-s_1\alpha_2}} 
&1\cdot s_2\cdot s_1 \ar[dl]^{y_{-s_1s_2\alpha_1}} \\
 &s_1\cdot s_2\cdot s_1
 }
&\xymatrix{ \\ \\ \mapright{(\gamma_{121})_!} \\ \\ }
&\xymatrix{
  &1 \ar[dl]_{y_{-\alpha_1}} \ar[dr]^{y_{-\alpha_2}} \ar[ddd]^<<<<<<<<<{y_{-(\alpha_1+\alpha_2)}} \\
s_1 \ar[d]_{y_{-\alpha_2}} \ar[drr]_<<<<<<<<<<{y_{-(\alpha_1+\alpha_2)}} & &s_2\ar[d]^{y_{-\alpha_1}}
\ar[dll]^<<<<<<<<<{y_{-(\alpha_1+\alpha_2)} } \\
s_1s_2 \ar[dr]_{y_{-\alpha_1}} & &s_2s_1 \ar[dl]^{y_{-\alpha_2}} \\
 &s_1s_2s_1=s_2s_1s_2
 }
\\ \\
\mapup{\tau^*} &&\mapup{\pi_{1}^*} 
\\ \\
  \xymatrix{
  &1\cdot 1 \ar[dl]_{y_{-\alpha_1}}  \ar[d]^>>{y_{-\alpha_2}} \\
s_1\cdot 1 \ar[d]_{y_{-\alpha_2}} 
&1\cdot s_2 \ar[dl]^<<<{y_{-s_2\alpha_1}} 
& \\ 
s_1\cdot s_2 
& 
& \\
 &
 }
&\xymatrix{ \\ \\ \mapright{(\gamma_{12})_!} \\ \\ } 
\xymatrix{
  &1 \ar[dl]_{y_{-\alpha_1}} \ar[dr]^{y_{-\alpha_2}} \ar[ddd]^<<<<<<<<<{y_{-(\alpha_1+\alpha_2)}} \\
s_1 \ar[d]_{y_{-\alpha_2}} \ar[drr]_<<<<<<<<<<{y_{-(\alpha_1+\alpha_2)}} & &s_2\ar[d]^{y_{-\alpha_1}}
\ar[dll]^<<<<<<<<<{y_{-(\alpha_1+\alpha_2)} } \\
s_1s_2 \ar[dr]_{y_{-\alpha_1}} & &s_2s_1 \ar[dl]^{y_{-\alpha_2}} \\
 &s_1s_2s_1=s_2s_1s_2
 }
\xymatrix{ \\ \\ \mapright{(\pi_1)_!} \\ \\ }
&\xymatrix{
  &1 \ar[dr]^{y_{-\alpha_2}} \ar[ddr]_{y_{-(\alpha_1+\alpha_2)}} \\ 
& &s_2 \ar[d]^{y_{-\alpha_1}}  \\
& & s_2s_1 \\
 &\hbox{$\phantom{ s_1s_2s_1=s_2s_1s_2 }$}
} 
\end{matrix}
$$
\caption{An example of the moment graphs for the diagram \eqref{BEexample}}
\label{BGGmomgraphs}
\end{sidewaysfigure}

\medskip\noindent
Theorem \ref{BSclassthm} says that the values on the vertices of the element $[Z_{i_1\cdots i_\ell}]$
on the moment graph of $\Gamma_{i_1,\ldots, i_\ell}$
are exactly the coefficients of the $2^\ell$ terms in the expansion of
$$A_{i_1}\cdots A_{i_\ell}
=(1+t_{s_{i_1}})\frac{1}{x_{-\alpha_{i_1}}}\cdots (1+t_{s_{i_\ell}})\frac{1}{x_{-\alpha_{i_\ell}}}.
$$
For example, in type $GL_3$,
$$[Z_{121}] = \left(
\begin{array}{c}
\frac{y_{-(\alpha_1+\alpha_2)}}{y_{-\alpha_1}}1\cdot1\cdot1 \\
+\frac{y_{-\alpha_2}}{y_{\alpha_1}}t_{s_1}\cdot 1\cdot 1
+1\cdot t_{s_2}\cdot 1 
+\frac{y_{-(\alpha_1+\alpha_2)}}{y_{-\alpha_1}}1\cdot1\cdot t_{s_1} \\
+t_{s_1}\cdot t_{s_2}\cdot 1 
+\frac{y_{-\alpha_2}}{y_{\alpha_1}}t_{s_1}\cdot 1\cdot t_{s_1} +1\cdot t_{s_2}\cdot t_{s_1} \\
+ t_{s_1}\cdot t_{s_2}\cdot t_{s_1}
\end{array}\right) b_1
$$
provides the expansion of 
$[Z_{121}] = (1+t_{s_1})\frac{1}{x_{-\alpha_1}}(1+t_{s_1})\frac{1}{x_{-\alpha_1}}(1+t_{s_1})\frac{1}{x_{-\alpha_1}}
y_{R^-}b_1$ in the basis $\{b_w\ |\ w\in W_0\}$.
An example of the pushpull in \eqref{BEdiagram} in the case of type $GL_3$
\begin{equation}\label{BEexample}
\xymatrix{P_1\times_B P_2\times_B P_1\times_B \pt
\ar^>>>>>>>>>>>>>>>{\gamma_{121}}[rr] \ar[d]_{\tau}&&GL_3/B \ar^{\pi_1}[d] \\
P_1 \times_B P_2\times_B\pt 
 \ar[r]_>>>>>>>{\gamma_{12}} &GL_3/B \ar_{\pi_1}[r]&GL_3/P_1
}
\end{equation}
has moment graphs as in Figure \ref{BGGmomgraphs},
and the computation in \eqref{BGGsquare} for this example is 
$$
\begin{matrix}
\begin{matrix}
&1 \\
1 &1 & 1 \\
1 &1 & 1 \\
&1
\end{matrix}
&\mapright{(\gamma_{121})_!}
&\begin{matrix}
&\Delta_{121} \\
\Delta_{121} &&1 \\
1 &&1 \\
&1
\end{matrix}
\\ \\
\mapup{\tau^*} &&\mapup{\pi_{1}^*} 
\\ \\
\begin{matrix}
&1 \\
1 &1 & \phantom{1} \\
1 &\phantom{1} & \phantom{1} \\
&\phantom{1}
\end{matrix}
&\mapright{(\gamma_{12})_!} \qquad
\begin{matrix}
&y_{-(\alpha_1+\alpha_2)} \\
y_{-\alpha_2} & &y_{-(\alpha_1+\alpha_2)} \\
0  &&y_{-\alpha_2} \\
&0
\end{matrix}
\qquad\mapright{(\pi_1)_!}
&\begin{matrix}
&\Delta_{121} \\
\phantom{?} &&1 \\
\phantom{?} &&1 \\
&\phantom{?}
\end{matrix}
\end{matrix}
$$
where $\Delta_{121} = \frac{y_{-(\alpha_1+\alpha_2)}}{y_{-\alpha_1}}+\frac{y_{-\alpha_2}}{y_{\alpha_1}}$.

\subsection{Change of groups morphisms across $\iota\colon B\hookrightarrow P_J$ }

In the same way that Theorem \ref{GBmom} provides $S\otimes_S^{W_0} S \cong \Omega_T(G/B)$ one can obtain
$$S^{W_J}\otimes_{S^W_0} S \cong \Omega_{P_J}(G/B),$$
and, if $\iota\colon B\hookrightarrow P_J$ is the inclusion
then the change of group homomorphisms 
$$
\iota^J\colon \Omega_{P_J}(G/B)\to \Omega_T(G/B)
\qquad\hbox{and}\qquad
\iota_J\colon \Omega_T(G/B)\to \Omega_{P_J}(G/B) 
$$
are given, combinatorially, by
$$\iota^J\colon S^{W_J}\otimes_{S^W_0} S \hookrightarrow S\otimes_{S^W_0} S
\qquad\hbox{and}\qquad
\iota_J\colon S\otimes_{S^W_0} S \longrightarrow S^{W_J}\otimes_{S^W_0} S,$$
with
$$\iota^J(f\otimes g) = \sum_{w\in W_J} w\left(\frac{1}{y_J}f\right)\otimes g,
\qquad\hbox{where}\quad
y_J = \prod_{\alpha\in R_J^+} y_{-\alpha},
$$
with $R_J^+$ the set of positive roots for $P_J\supseteq B\supseteq T$.  The pushforward
$\iota^J$ is similar to the pushforward operator $(\pi_J)_!$ appearing in \eqref{GPpushpul[}
except acting on the left factor of $S\otimes_{S^{W_0}} S$ (see, for example, the definition 
of $\delta_i$ in \cite[\S 7]{Ka}).

\section{Schubert classes $[X_w]$}\label{Schubertclasses}

Now we consider the inclusions 
$\sigma_w\colon X_w \longrightarrow G/B$
of the Schubert varieties into the flag variety.
For $w\in W_0$, define the \emph{Schubert classes}\quad
\begin{equation}\label{Schubdefn}
[X_w] = (\sigma_w)_!(1),
\qquad\hbox{where}\quad
(\sigma_w)_!\colon \Omega_T(X_w)\to \Omega_T(G/B).
\end{equation}
If $X_w$ is not smooth then, as discussed further below, it is not clear that $(\sigma_w)_!$ is well defined.
Though we consider various approaches to the analysis of $[X_w] = (\sigma_w)_!(1)$ below, we have not yet found a definition
of $(\sigma_w)_!$ which is fully satisfying (at least to us) in the singular case.

In \emph{generalized} cohomology 
$$\hbox{the Schubert class $[X_w]$ is \emph{not always} equal to $[Z_{\vec w}]$}
$$
for a reduced word $\vec w$ of $w$, although, in equivariant cohomology and 
equivariant $K$-theory, $[X_w]=[Z_{\vec w}]$ if 
$\vec w$ is a reduced word for $w$.  We consider various approaches to the analysis of $[X_w] = (\sigma_w)_!(1)$:
\begin{enumerate}
\item[(a)] Defining $(\sigma_w)_!(1)$ by \eqref{pullandpush};
\item[(b)] Comparing $[X_w] = (\sigma_w)_!(1)$ and the Bott-Samelson class $[Z_{\vec w}]$ via the diagram
\begin{equation}\label{compXandZ}
\xymatrix{
\Omega_T(\Gamma_{\vec w}) \ar[dr]^{(\gamma_{\vec w})_!} \ar[d]_{(\tilde\gamma_{\vec w})_!} \\
\Omega_T(X_w) \ar[r]^{(\sigma_w)_!} &\Omega_T(G/B)
}
\end{equation}
\item[(c)] Combinatorial forcing by support conditions, normalization and/or $(S,S)$-bimodule structure of the cohomology.
\end{enumerate}

\noindent
\textbf{(a) Is $(\sigma_w)_!(1)$ given by \eqref{pullandpush}?} As pointed out in \cite[Proposition 2.7]{Ty}, since $X_w$ is filtered
by Schubert cells $BvB$ with $v\le w$ and $BvB\cong \CC^{\ell(v)}$ has even
real dimension, the Schubert variety $X_w$ has no odd-dimensional cohomology, and thus, by
\cite[Theorem 14]{GKM}, the Schubert variety $X_w$ is `equivariantly formal' (i.e., is a GKM-space)
and the moment graph theory applies.  The moment graph of $X_w$ is the subgraph of the moment graph of $G/B$
with vertices $\{ v\in W_0\ |\ v\le w\}$.   If $X_w$ is smooth then there are no challenges in
defining the pushforward $(\sigma_w)_!$ and the pushforward formula
in \eqref{pullandpush} gives that 
\begin{equation}\label{smSchclass}
\hbox{if $X_w$ is smooth,}\qquad
\hbox{then}\quad [X_w]_v = \frac{y_{R^-}}{\displaystyle{
\prod_{\beta\in R^+\atop vs_\beta\le w} y_{-\beta} }},
\quad\hbox{for $v\in W_0$ such that $v\le w$,}
\end{equation}
as found, for example, in \cite[Theorem 7.2.1]{BiLa} (the notation $f=\sum_{w\in W_0} f_wb_w$ for elements of $\Omega_T(G/B)$
is as \eqref{momgrnotation}). 
For example, the inclusion $\sigma_{s_2s_1}\colon X_{s_2s_1}\to G/B$ for $G=GL_3$
corresponds to the inclusion of moment graphs
$$
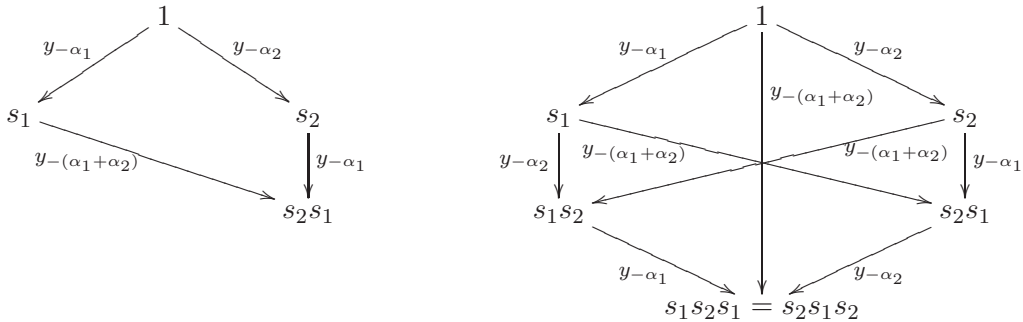
\xymatrix{
&1 \ar[dl]_{y_{-\alpha_1}} \ar[dr]^{y_{-\alpha_2}} \\ 
s_1
 \ar[drr]_<<<<<<<<<<{y_{-(\alpha_1+\alpha_2)}} & &s_2 \ar[d]^{y_{-\alpha_1}} \\
\hbox{$\phantom{s_1s_2}$} 
& &s_2s_1  \\ 
&\hbox{$\phantom{s_1s_2s_1}$}
}
\qquad\qquad
\xymatrix{
  &1 \ar[dl]_{y_{-\alpha_1}} \ar[dr]^{y_{-\alpha_2}} \ar[ddd]^<<<<<<<<<{y_{-(\alpha_1+\alpha_2)}} \\
s_1 \ar[d]_{y_{-\alpha_2}} \ar[drr]_<<<<<<<<<<{y_{-(\alpha_1+\alpha_2)}} & &s_2\ar[d]^{y_{-\alpha_1}}
\ar[dll]^<<<<<<<<<{y_{-(\alpha_1+\alpha_2)} } \\
s_1s_2 \ar[dr]_{y_{-\alpha_1}} & &s_2s_1 \ar[dl]^{y_{-\alpha_2}} \\
 &s_1s_2s_1=s_2s_1s_2
}
$$
so that
$$[X_{s_2s_1}]\quad=\quad
\begin{matrix}
&\frac{y_{R^-}}{y_{-\alpha_1}y_{-\alpha_2}} \\
\frac{y_{R^-}}{y_{-\alpha_1}y_{-s_1\alpha_2}} &&\frac{y_{R^-}}{y_{-\alpha_1}y_{-\alpha_2}} \\
\\
0 &&\frac{y_{R^-}}{y_{-\alpha_1}y_{-s_1\alpha_2}} \\
 &0 & 
\end{matrix}
$$

The following example illustrates that this procedure does not work well when
$X_w$ is not smooth. From \cite[Prop. 6.1]{Ku}, the singular Schubert varieties for $G$ of rank 2 are
$$\begin{array}{ccccc}
\hbox{Type} &\phantom{TTTTTTT} &\hbox{Singular} &\phantom{TTTTTTT} &\hbox{Locus} \\
\\
B_2 &&X_{s_1s_2s_1} &&X_{s_1} \\
\\
G_2 &&X_{s_1s_2s_1} &&X_{s_1} \\
G_2 &&X_{s_1s_2s_1s_2} &&X_{s_1s_2} \\
G_2 &&X_{s_2s_1s_2s_1} &&X_{s_2s_1} \\
G_2 &&X_{s_1s_2s_1s_2s_1} &&X_{s_1s_2s_1} \\
G_2 &&X_{s_2s_1s_2s_1s_2} &&X_{s_2} 
\end{array} 
$$
The inclusion $\sigma_{s_1s_2s_1}\colon X_{s_1s_2s_1}\to G/B$ for $G=Sp_4$ (Type $B_2$) corresponds to the inclusion of
moment graphs
$$
\xymatrix{
  &1 \ar[dl]_{y_{-\alpha_1}} \ar[dr]^{y_{-\alpha_2}} 
\ar[dddl]^<<<<<<<{y_{-s_1\alpha_2}} 
  \\
s_1 \ar[d]_{y_{-\alpha_2}} 
\ar[drr]^<<<<<<<<<<<<<<{y_{-s_1\alpha_2}} 
& &s_2\ar[d]^{y_{-\alpha_1}}
\ar[dll]_<<<<<<<<<<<<<<{y_{-s_2\alpha_1} } \\
s_1s_2 \ar[d]_{y_{-\alpha_1}}  
&  &s_2s_1 \ar[dll]_<<<<<<<<<<<<{y_{-s_2\alpha_1}} \\ 
s_1s_2s_1 
& &\hbox{$\phantom{s_2s_1s_2}$} \\ 
 &\hbox{$\phantom{s_1s_2s_1s_2}$}
}
\qquad\qquad
\xymatrix{
  &1 \ar[dl]_{y_{-\alpha_1}} \ar[dr]^{y_{-\alpha_2}} 
\ar[dddl]_<<<<<<<<<<<{y_{-s_1\alpha_2}} 
\ar[dddr]^<<<<<<<<<<<{y_{-s_2\alpha_1}} 
  \\
s_1 \ar[d]_{y_{-\alpha_2}} 
\ar[dddr]^>>>>>>>>>>>>{y_{-s_2\alpha_1}} 
\ar[drr]^<<<<<<<<<<<<<<{y_{-s_1\alpha_2}} 
& &s_2\ar[d]^{y_{-\alpha_1}}
\ar[dddl]^>>>>>>>>>{y_{-s_1\alpha_2} } 
\ar[dll]_<<<<<<<<<<<<<<{y_{-s_2\alpha_1} } \\
s_1s_2 \ar[drr]^<<<<<<<<<<<<{y_{-s_1\alpha_2}} \ar[d]_{y_{-\alpha_1}} & 
&s_2s_1 \ar[dll]_<<<<<<<<<<<<{y_{-s_2\alpha_1}} \ar[d]^{y_{-\alpha_2}} \\
s_1s_2s_1 \ar[dr]_{y_{-\alpha_2}} & &s_2s_1s_2 \ar[dl]^{y_{-\alpha_1}} \\
 &s_1s_2s_1s_2
}
$$
but the direct ``naive'' application of the pushforward formula in \eqref{pullandpush} produces
\begin{equation}\label{OOPS}
[X_{s_1s_2s_1}]?=?\qquad
\begin{matrix}
&y_{-s_2\alpha_1} \\
y_{-s_2\alpha_1} &&y_{-s_1\alpha_2} \\
y_{-s_1\alpha_2} &&y_{-\alpha_2} \\
y_{-\alpha_2} & &0 \\
 &0 & 
\end{matrix}
~=~
\begin{matrix}
&y_{-(\alpha_1+\alpha_2)} \\
y_{-(\alpha_1+\alpha_2)} &&y_{-(2\alpha_1+\alpha_2)} \\
y_{-(2\alpha_1+\alpha_2)} &&y_{-\alpha_2} \\
y_{-\alpha_2} & &0 \\
 &0 & 
\end{matrix}
\end{equation}
which cannot be correct for $[X_{s_1s_2s_1}]$ since the right hand side
does not satisfy the condition to be in $\im\, \Phi$ (the difference across the edge $1\to s_2$ is not divisible by $y_{-\alpha_2}$).  This
answer needs to be corrected by finding $N$ so that
$$[X_{s_1s_2s_1}]
~=~
\begin{matrix}
&Ny_{-(\alpha_1+\alpha_2)} \\
Ny_{-(\alpha_1+\alpha_2)} &&y_{-(2\alpha_1+\alpha_2)} \\
y_{-(2\alpha_1+\alpha_2)} &&y_{-\alpha_2} \\
y_{-\alpha_2} & &0 \\
 &0 & 
\end{matrix}
$$
where the correction factor $N$ appears on vertices corresponding to the singular locus.

In the example in \eqref{OOPS} we see that the moment  graph knows that $X_{s_1s_2s_1}$ is not smooth!  
It is interesting to contrast \eqref{OOPS}
with the same analysis for $\sigma_{s_2s_1s_2}\colon X_{s_2s_1s_2}\to G/B$,
where the pushforward formula gives
$$[X_{s_2s_1s_2}]=\qquad
\begin{matrix}
&y_{-s_1\alpha_2} \\
y_{-s_2\alpha_1} &&y_{-s_1\alpha_2} \\
y_{-\alpha_1} &&y_{-s_2\alpha_1} \\
0 & &y_{-\alpha_1} \\
 &0 & 
\end{matrix}
~=~
\begin{matrix}
&y_{-(2\alpha_1+\alpha_2)} \\
y_{-(\alpha_1+\alpha_2)} &&y_{-(2\alpha_1+\alpha_2)} \\
y_{-\alpha_1} &&y_{-(\alpha_1+\alpha_2)} \\
0 & &y_{-\alpha_1} \\
 &0 & 
\end{matrix}
$$
which is in $\im\, \Phi$ (this case works out well since $X_{s_2s_1s_2}$ is smooth).

\medskip\noindent
\textbf{(b) Using \eqref{compXandZ} to compare $[X_w]$ and $[Z_{\vec w}]$.}
Working in rank 2, use notations $y_R^-$, $\Delta_{121}$ and $\Delta_{212}$ as in 
\eqref{Deltadefn},
so that (see \eqref{BEexample} and Figure \ref{BGGmomgraphs})
$$
[Z_{212}] \quad=\quad
\begin{matrix}
&\Delta_{212} \\
\frac{y_{R^-}}{y_{-\alpha_1}y_{-\alpha_2}y_{-s_1\alpha_2}} &&\Delta_{212} \\
\frac{y_{R^-}}{y_{-\alpha_2}y_{-s_2\alpha_1}y_{-s_2s_1\alpha_2}} &&\frac{y_{R^-}}{y_{-\alpha_1}y_{-\alpha_2}y_{-s_1\alpha_2}} \\
0 &&\frac{y_{R^-}}{y_{-\alpha_2}y_{-s_2\alpha_1}y_{-s_2s_1\alpha_2}} \\
 &0& \\
\end{matrix}
$$
Since $X_{s_1s_2s_1}$ is smooth it is reasonable to apply the pushforward formula
in \eqref{pullandpush} which gives
$$
[X_{s_2s_1s_2}]
\quad=\quad
\begin{matrix}
&\frac{y_{R^-}}{y_{-\alpha_1}y_{-\alpha_2}y_{-s_2\alpha_1}} \\
\frac{y_{R^-}}{y_{-\alpha_1}y_{-\alpha_2}y_{-s_1\alpha_2}} &&\frac{y_{R^-}}{y_{-\alpha_1}y_{-\alpha_2}y_{-s_2\alpha_1}} \\
\frac{y_{R^-}}{y_{-\alpha_2}y_{-s_2\alpha_1}y_{-s_2s_1\alpha_2}} &&\frac{y_{R^-}}{y_{-\alpha_1}y_{-\alpha_2}y_{-s_1\alpha_2}} \\
0 &&\frac{y_{R^-}}{y_{-\alpha_2}y_{-s_2\alpha_1}y_{-s_2s_1\alpha_2}} \\
 &0& \\
\end{matrix}
$$
Using these and computing with the formulas \eqref{yalphynegalapha}-\eqref{divdiffexpansion} gives the formula
\begin{align*}
[Z_{212}] 
&= [X_{s_2s_1s_2}] + \left(\Delta_{212}-\frac{y_{R^-}}{y_{-\alpha_1}y_{-\alpha_2}y_{-s_2\alpha_1}}\right)
\frac{y_{-\alpha_2}}{y_{R^-}}[X_{s_2}] \\
&= [X_{s_2s_1s_2}] + \frac{y_{R^-}}{y_{-\alpha_1}y_{-\alpha_2}y_{-s_2\alpha_1}}
\left( \frac{y_{-s_2\alpha_1}-y_{-\alpha_1}}{y_{-\alpha_2}}
+p(y_{\alpha_2},y_{-\alpha_2})y_{-\alpha_1}-1\right)\frac{y_{-\alpha_2}}{y_{R^-}}[X_{s_2}] \\
&= [X_{s_2s_1s_2}] + \frac{y_{R^-}}{y_{-\alpha_1}y_{-\alpha_2}y_{-s_2\alpha_1}}
\left( (1-p(y_{-\alpha_1},y_{-\alpha_2})y_{-\alpha_1} 
+p(y_{\alpha_2},y_{-\alpha_2})y_{-\alpha_1}-1\right)\frac{y_{-\alpha_2}}{y_{R^-}}[X_{s_2}] \\
&= [X_{s_2s_1s_2}] + \frac{1}{y_{-s_2\alpha_1}}
\big(p(y_{\alpha_2},y_{-\alpha_2})-p(y_{-\alpha_1},y_{-\alpha_2})\big)y_{-\alpha_1}\frac{y_{-\alpha_2}}{y_{R^-}}[X_{s_2}]
\end{align*}
which is reflected in \cite[17.3, first equation]{CPZ} and \cite[\S5.2]{HK}.
Similarly, with our conjectured correction factor $N$ as in \eqref{Ndefn}, we get a formula which would provide
$[Z_{121}]-[X_{s_1s_2s_1}]=0$ in cohomology and K-theory
but have $[Z_{121}]-[X_{s_1s_2s_1}]\ne 0$ in complex or algebraic cobordism:
\begin{align*}
[Z_{121}]&-[X_{s_1s_2s_1}] 
= \left(\Delta_{121}-\frac{Ny_{R^-}}{y_{-\alpha_1}y_{-\alpha_2}y_{-s_1\alpha_2}}\right)\frac{y_{-\alpha_1}}{y_{R^-}}[X_{s_1}] \\
&= \frac{y_{R^-}}{y_{-\alpha_1}y_{-\alpha_2}y_{-s_1\alpha_2}}
\left( \frac{y_{-s_1\alpha_2}-y_{-\alpha_2}}{y_{-\alpha_1}}
+p(y_{\alpha_1},y_{-\alpha_1})y_{-\alpha_2} - N \right)\frac{y_{-\alpha_1}}{y_{R^-}}[X_{s_1}] \\
&= \frac{y_{R^-}}{y_{-\alpha_1}y_{-\alpha_2}y_{-s_1\alpha_2}}
\left(
\begin{array}{l}
(1-p(y_{-\alpha_2},y_{-j\alpha_1})y_{-\alpha_2})\big(1+\sum_{k=1}^{j-1} (1-p(y_{-\alpha_1},y_{-k\alpha_1})y_{-k\alpha_1}\big) 
\\
\quad
+p(y_{\alpha_1},y_{-\alpha_1})y_{-\alpha_2} - N
\end{array}\right)\frac{y_{-\alpha_1}}{y_{R^-}}[X_{s_1}] \\
&=\frac{y_{R^-}}{y_{-\alpha_1}y_{-\alpha_2}y_{-s_1\alpha_2}}
\big(p(y_{\alpha_1},y_{-\alpha_1})-p(y_{-\alpha_2},y_{-j\alpha_1})\big)y_{-\alpha_2}
\frac{y_{-\alpha_1}}{y_{R^-}}[X_{s_1}] \\
&=\frac{1}{y_{-s_1\alpha_2}}
\big(p(y_{\alpha_1},y_{-\alpha_1})-p(y_{-\alpha_2},y_{-j\alpha_1})\big)[X_{s_1}].
\end{align*}

\smallskip\noindent
\textbf{(c) Combinatorial forcing:} The Schubert classes satisfy
\begin{enumerate}
\item[(a)] (normalization) $[X_w]_w=\prod_{\alpha\in R(w)} y_{-\alpha}$, where $R(w) = \{ \alpha\in R^+\ |\ w\alpha\not\in R^+\}$.
\item[(b)] If $W_Ju=W_Jz$ then $[X_{w_Ju}]_v=[X_{w_Ju}]_z$,
\item[(c)] (support) $[X_w]_v=0$ unless $v\le w$.
\end{enumerate}
These properties do not characterize the Schubert classes; the Bott-Samelson classes also satisfy these
properties.  
As observed, for example, in \cite[Proposition 4.3]{HHH}, in equivariant cohomology
a degree condition can be imposed to get uniqueness.
It is not clear to us how to generalize the degree condition to equivariant K-theory and/or equivariant
cobordism. It seems plausible that in generalized equivariant cohomology the Schubert classes
might be characterized by positivity properties,
or by using the $(S,S)$-bimodule structure of $\Omega_T(X_w)$ and $\Omega_T(Z_{\vec w})$
as in the theory of Soergel bimodules (see \cite{Soe} and \cite{EW}).

\section{Products with Schubert classes}\label{Schubertproducts}

For $w\in W_0$ define Schubert classes $[X_w]$ by $[X_w]=(\sigma_w)_!(1)$ as in \eqref{Schubdefn}.  Continue to use
notations $f = \sum_{w\in W_0} f_w b_w$ for elements of $\Omega_T(G/B)$, as in \eqref{momgrnotation}.

\medskip\noindent
{\bf The Schubert product problem}:  {\sl Find a combinatorial description of the 
$c_{uv}^w\in R$ given by 
\begin{equation}\label{Schubproblem}
[X_u][X_v] = \sum_{w\in W_0} c_{uv}^w [X_w].
\end{equation}
}
As is visible from the formula \eqref{Schubexpansion} below and the formulas at the end of this section,
if $v\le u$ in Bruhat order then
\begin{equation}\label{refinedSchubproblem}
[X_u][X_v] = [X_u]_v [X_v] + \sum_{w<v} c_{uv}^w [X_w],
\end{equation}
and so the determination of the moment graph values $[X_u]_v$ is a subproblem of the Schubert product problem.
The other coefficients $c_{uv}^w$ are determined by the $[X_u]_v$ in an intricate but, perhaps, controllable fashion.
Furthermore, our computations of products in the rank two cases display a certain amount of positivity, indicating that 
there may be a positivity statement for equivariant cobordism analogous to that which holds for equivariant
cohomology and equivariant K-theory (see \cite{Gra} and \cite{AGM}).

Properties (a) and (c) are already enough to provide an algorithm for expanding
an element $f = \sum_{w\in W_0} f_w b_w$ in terms of Schubert classes.
If $f$ has support on $w$ with $\ell(w)\le k$ then
$$f - \sum_{\ell(w)=k} f_w\frac{1}{[X_w]_w}[X_w]
= \sum_{\ell(v)\le k-1}\left(f_v-\sum_{\ell(w)=k} f_w\frac{[X_w]_v}{[X_w]_w}\right) b_v
$$
has support on $v$ with $\ell(v)\le k-1$.  Then
\begin{align*}
f - &\sum_{\ell(w)=k} f_w\frac{1}{[X_w]_w}[X_w]
-\sum_{\ell(v)=k-1} \left(f_v
- \sum_{\ell(v)=k-1\atop \ell(w)=k} f_w\frac{[X_w]_v}{[X_w]_w}\right)\frac{1}{[X_v]_v}[X_v]
\\
&= \sum_{\ell(z)\le k-2}
\left(
f_z
-\sum_{\ell(w)=k} f_w\frac{[X_w]_z}{[X_w]_w} 
-\sum_{\ell(v)=k-1} f_v\frac{[X_v]_z}{[X_v]_v} 
+ \sum_{\ell(v)=k-1\atop \ell(w)=k} f_w\frac{[X_w]_v}{[X_w]_w}\frac{[X_v]_z}{[X_v]_v}
\right)b_z
\end{align*}
and induction gives that
\begin{equation}\label{Schubexpansion}
f = \sum_{z\in W_0}
\left( \sum_{k=1}^{\ell(w_0)} \sum_{w_1>\cdots> w_k=z} (-1)^{k-1} f_{w_1}
\frac{[X_{w_1}]_{w_2}}{[X_{w_1}]_{w_1}}
\frac{[X_{w_2}]_{w_3}}{[X_{w_2}]_{w_2}}
\cdots
\frac{[X_{w_{k-1}}]_{w_k}}{[X_{w_{k-1}}]_{w_{k-1}}}
\frac{1}{[X_{w_k}]_{w_k}}\right) [X_z]
\end{equation}
with the terms in the sum naturally indexed by chains in the
Bruhat order (compare to, for example, \cite{BS}).

For example, in rank 2 using notations as in Section \ref{rank2comps}, if $f = \sum_{w\le s_1s_2s_1s_2} f_w b_w$ then
\begin{align*}
f &= f_{s_1s_2s_1s_2}\frac{1}{[X_{s_1s_2s_1s_2}]_{s_1s_2s_1s_2}}[X_{s_1s_2s_1s_2}] \\
&\quad
+(f_{s_1s_2s_1}-f_{s_1s_2s_1s_2})\frac{1}{[X_{s_1s_2s_1}]_{s_1s_2s_1}} [X_{s_1s_2s_1}]
+(f_{s_2s_1s_2}-f_{s_1s_2s_1s_2})\frac{1}{[X_{s_2s_1s_2}]_{s_2s_1s_2}} [X_{s_2s_1s_2}] \\
&\quad+
\left(
(f_{s_1s_2}-f_{s_2s_1s_2})
+(f_{s_1s_2s_1s_2}-f_{s_1s_2s_1})\frac{[X_{s_1s_2s_1}]_{s_1s_2}}{[X_{s_1s_2s_1}]_{s_1s_2s_1}}
\right)\frac{1}{[X_{s_1s_2}]_{s_1s_2}} [X_{s_1s_2}] \\
&\quad+
\left(
(f_{s_2s_1}-f_{s_1s_2s_1})
+(f_{s_1s_2s_1s_2}-f_{s_2s_1s_2})\frac{[X_{s_2s_1s_2}]_{s_2s_1}}{[X_{s_2s_1s_2}]_{s_2s_1s_2}}
\right)\frac{1}{[X_{s_2s_1}]_{s_2s_1}} [X_{s_2s_1}] \\
&\quad+
\left(\begin{array}{l}
(f_{s_1}-f_{s_2s_1})+(f_{s_2s_1s_2}-f_{s_1s_2})\frac{[X_{s_1s_2}]_{s_1}}{[X_{s_1s_2}]_{s_1s_2}} \\
\quad+(f_{s_1s_2s_1s_2}-f_{s_1s_2s_1})\bigg(
\frac{[X_{s_1s_2s_1}]_{s_1}}{[X_{s_1s_2s_1}]_{s_1s_2s_1}}
-\frac{[X_{s_1s_2s_1}]_{s_1s_2}}{[X_{s_1s_2s_1}]_{s_1s_2s_1}}
\frac{[X_{s_1s_2}]_{s_1}}{[X_{s_1s_2}]_{s_1s_2}}
-1\bigg)
\end{array}
\right)\frac{1}{[X_{s_1}]_{s_1}} [X_{s_1}] \\
&\quad+
\left(\begin{array}{l}
(f_{s_2}-f_{s_1s_2})+(f_{s_1s_2s_1}-f_{s_2s_1})\frac{[X_{s_2s_1}]_{s_2}}{[X_{s_2s_1}]_{s_2s_1}} \\
\quad
+(f_{s_2s_1s_2s_1}-f_{s_2s_1s_2})\bigg(
\frac{[X_{s_2s_1s_2}]_{s_2}}{[X_{s_2s_1s_2}]_{s_2s_1s_2}}
-\frac{[X_{s_2s_1s_2}]_{s_2s_1}}{[X_{s_2s_1s_2}]_{s_2s_1s_2}}
\frac{[X_{s_2s_1}]_{s_2}}{[X_{s_2s_1}]_{s_2s_1}}
-1\bigg)
\end{array}
\right)\frac{1}{[X_{s_2}]_{s_2}} [X_{s_2}] \\
&\quad+
(f_1-f_{s_1}-f_{s_2}+f_{s_1s_2}+f_{s_2s_1}-f_{s_1s_2s_1}-f_{s_2s_1s_2}+f_{s_1s_2s_1s_2})
\frac{1}{[X_1]_1}[X_1]
\end{align*}
and we may use the explicit values of $[X_w]_v$ given in Figure \ref{rank2Schuberts} to derive
\begin{align*}
f &= f_{s_1s_2s_1s_2}\frac{y_{-\alpha_1}y_{-s_1\alpha_2}y_{-s_1s_2\alpha_1}y_{-s_1s_1s_1\alpha_2}}{y_{R^-}}[X_{s_1s_2s_1s_2}] \\
&\quad
+(f_{s_1s_2s_1}-f_{s_1s_2s_1s_2})
\frac{y_{-\alpha_1}y_{-s_1\alpha_2}y_{-s_1s_2\alpha_1}}{y_{R^-}} [X_{s_1s_2s_1}]
+(f_{s_2s_1s_2}-f_{s_1s_2s_1s_2})
\frac{y_{-\alpha_2}y_{-s_2\alpha_1}y_{-s_2s_1\alpha_2}}{y_{R^-}} [X_{s_2s_1s_2}] \\
&\quad+
\left(
(f_{s_1s_2}-f_{s_2s_1s_2})
+(f_{s_1s_2s_1s_2}-f_{s_1s_2s_1})
\frac{y_{-\alpha_1}y_{-s_1\alpha_2}y_{-s_1s_2\alpha_1}}{y_{-\alpha_1}y_{-\alpha_2}y_{-s_2\alpha_1}}
\right)\frac{y_{-\alpha_2}y_{-s_2\alpha_1}}{y_{R^-}} [X_{s_1s_2}] \\
&\quad+
\left(
(f_{s_2s_1}-f_{s_1s_2s_1})
+(f_{s_1s_2s_1s_2}-f_{s_2s_1s_2})
\frac{y_{-\alpha_2}y_{-s_2\alpha_1}y_{-s_2s_1\alpha_2}}{y_{-\alpha_1}y_{-\alpha_2}y_{-s_1\alpha_2}}
\right)\frac{y_{-\alpha_1}y_{-s_1\alpha_2}}{y_{R^-}} [X_{s_2s_1}] \\
&\quad+
\left(\begin{array}{l}
(f_{s_1}-f_{s_2s_1})+(f_{s_2s_1s_2}-f_{s_1s_2})
\frac{y_{-\alpha_2}y_{-s_2\alpha_1}}{y_{-\alpha_1}y_{-\alpha_2}} \\
\quad+(f_{s_1s_2s_1s_2}-f_{s_1s_2s_1})\bigg(
\frac{Ny_{-\alpha_1}y_{-s_1\alpha_2}y_{-s_1s_2\alpha_1}}{y_{-\alpha_1}y_{-\alpha_2}y_{-s_1\alpha_2}}
-\frac{y_{-\alpha_1}y_{-s_1\alpha_2}y_{-s_1s_2\alpha_1}}{y_{-\alpha_1}y_{-\alpha_2}y_{-s_2\alpha_1}}
\frac{y_{-\alpha_2}y_{-s_2\alpha_1}}{y_{-\alpha_1}y_{-\alpha_2}}
-1\bigg)
\end{array}
\right)\frac{y_{-\alpha_1}}{y_{R^-}} [X_{s_1}] \\
&\quad+
\left(\begin{array}{l}
(f_{s_2}-f_{s_1s_2})+(f_{s_1s_2s_1}-f_{s_2s_1})
\frac{y_{-\alpha_1}y_{-s_1\alpha_2}}{y_{-\alpha_1}y_{-\alpha_2}} \\
\quad
+(f_{s_2s_1s_2s_1}-f_{s_2s_1s_2})\bigg(
\frac{y_{-\alpha_2}y_{-s_2\alpha_1}y_{-s_2s_1\alpha_2}}{y_{-\alpha_1}y_{-\alpha_2}y_{-s_2\alpha_1}}
-\frac{y_{-\alpha_2}y_{-s_2\alpha_1}y_{-s_2s_1\alpha_2}}{y_{-\alpha_1}y_{-\alpha_2}y_{-s_1\alpha_2}}
\frac{y_{-\alpha_1}y_{-s_1\alpha_2}}{y_{-\alpha_1}y_{-\alpha_2}}
-1\bigg)
\end{array}
\right)\frac{y_{-\alpha_2}}{y_{R^-}} [X_{s_2}] \\
&\quad+
(f_1-f_{s_1}-f_{s_2}+f_{s_1s_2}+f_{s_2s_1}-f_{s_1s_2s_1}-f_{s_2s_1s_2}+f_{s_1s_2s_1s_2})
\frac{1}{y_{R^-}}[X_1]
\end{align*}
which simplifies to
\begin{align*}
y_{R^-}f &= f_{s_1s_2s_1s_2}y_{-\alpha_1}y_{-s_1\alpha_2}y_{-s_1s_2\alpha_1}y_{-s_1s_1s_1\alpha_2}[X_{s_1s_2s_1s_2}] \\
&\quad
+(f_{s_1s_2s_1}-f_{s_1s_2s_1s_2})
y_{-\alpha_1}y_{-s_1\alpha_2}y_{-s_1s_2\alpha_1} [X_{s_1s_2s_1}] \\
&\quad+(f_{s_2s_1s_2}-f_{s_1s_2s_1s_2})
y_{-\alpha_2}y_{-s_2\alpha_1}y_{-s_2s_1\alpha_2} [X_{s_2s_1s_2}] \\
&\quad+
(
(f_{s_1s_2}-f_{s_2s_1s_2})y_{-\alpha_2}y_{-s_2\alpha_1}
+(f_{s_1s_2s_1s_2}-f_{s_1s_2s_1})
y_{-s_1\alpha_2}y_{-s_1s_2\alpha_1}) [X_{s_1s_2}] \\
&\quad+
(
(f_{s_2s_1}-f_{s_1s_2s_1})y_{-\alpha_1}y_{-s_1\alpha_2}
+(f_{s_1s_2s_1s_2}-f_{s_2s_1s_2})
y_{-s_2\alpha_1}y_{-s_2s_1\alpha_2}) [X_{s_2s_1}] \\
&\quad+
\left(\begin{array}{l}
(f_{s_1}-f_{s_2s_1})y_{-\alpha_1}+(f_{s_2s_1s_2}-f_{s_1s_2})
y_{-s_2\alpha_1} \\
\quad+(f_{s_1s_2s_1s_2}-f_{s_1s_2s_1})\bigg(
\frac{Ny_{-s_1s_2\alpha_1}y_{-\alpha_1}}{y_{-\alpha_2}}
-\frac{y_{-s_1\alpha_2}y_{-s_1s_2\alpha_1}}{y_{-\alpha_2}}
-y_{-\alpha_1}\bigg)
\end{array}
\right) [X_{s_1}] \\
&\quad+
\left(\begin{array}{l}
(f_{s_2}-f_{s_1s_2})y_{-\alpha_2}+(f_{s_1s_2s_1}-f_{s_2s_1})
y_{-s_1\alpha_2} \\
\quad
+(f_{s_2s_1s_2s_1}-f_{s_2s_1s_2})\bigg(
\frac{y_{-s_2s_1\alpha_2}y_{-\alpha_2}}{y_{-\alpha_1}}
-\frac{y_{-s_2\alpha_1}y_{-s_2s_1\alpha_2}}{y_{-\alpha_1}}
-y_{-\alpha_2}\bigg)
\end{array}
\right) [X_{s_2}] \\
&\quad+
(f_1-f_{s_1}-f_{s_2}+f_{s_1s_2}+f_{s_2s_1}-f_{s_1s_2s_1}-f_{s_2s_1s_2}+f_{s_1s_2s_1s_2})
[X_1]
\end{align*}

This last formula allows for quick computation of products with Schubert classes in rank 2
for low dimensional Schubert varieties.  In particular,
for $g = \sum_{w\in W_0} g_w b_w$ in $\Omega_T(G/B)$, 
\begin{align*}
g[X_1] &= g_1[X_1], \\
g[X_{s_1}] &= g_{s_1}[X_{s_1}] + g_{1,s_1}[X_1], \qquad
\hbox{where}\quad
g_{1,s_1} = \frac{g_1-g_{s_1}}{y_{-\alpha_1}}, \\
g[X_{s_2}] &= g_{s_2}[X_{s_2}] +g_{1,s_2}[X_1], \qquad
\hbox{where}\quad
g_{1,s_2} = \frac{g_1-g_{s_2}}{y_{-\alpha_2}}, \\
g[X_{s_1s_2}] &=g_{s_1s_2}[X_{s_1s_2}] + g_{s_1,s_1s_2}[X_{s_1}]+g_{s_2,s_1s_2}[X_{s_2}]
+\frac{g_{1,s_1}-g_{s_2,s_1s_2}}{y_{-\alpha_2}}[X_1], \\
g[X_{s_2s_1}] &= g_{s_2s_1}[X_{s_2s_1}] + g_{s_1,s_2s_1} [X_{s_1}] + g_{s_2,s_2s_1}[X_{s_2}]
+\frac{g_{1,s_2}-g_{s_1,s_2s_1}}{y_{-\alpha_1}}[X_1],
\end{align*}
where
$$g_{s_1,s_1s_2} = \frac{g_{s_1}-g_{s_1s_2}}{y_{-\alpha_2}}, \quad
g_{s_2,s_1s_2} = \frac{g_{s_2}-g_{s_1s_2}}{y_{-s_2\alpha_1}}, \quad
g_{s_1,s_2s_1} = \frac{g_{s_1}-g_{s_2s_1}}{y_{-s_1\alpha_2}}, \quad
g_{s_2,s_2s_1} = \frac{g_{s_2}-g_{s_2s_1}}{y_{-\alpha_1}}.
$$
Using \eqref{nilhmomgr1}, Pieri-Chevalley rules giving the expansions of products $x_\lambda [X_w]$ in terms of
Schubert classes are directly determined from these formulas.


\section{Schubert classes and products in rank 2}\label{rank2comps}

In rank 2, $W_0$ is a dihedral group generated by $s_1$ and $s_2$ with $s_i^2=1$, $s_1\alpha_1=-\alpha_1$,
$s_2\alpha_2=-\alpha_2$, 
$$
\begin{matrix}
\begin{array}{ll}
s_1\alpha_1=-\alpha_1, &s_1\alpha_2 = j\alpha_1+\alpha_2, \\
s_2\alpha_1 = \alpha_1+\alpha_2, &s_2\alpha_2=-\alpha_2, \\
\end{array}
\\ \\
\hbox{with} \\ \\
j = \begin{cases} 
1, \hbox{in Type $A_2$,} \\
2, \hbox{in Type $B_2$,} \\
3, \hbox{in Type $G_2$,}
\end{cases}
\end{matrix}
\qquad\hbox{and}\qquad
\begin{matrix}
\begin{matrix}
&b_1 \\
b_{s_1} &&b_{s_2} \\
b_{s_1s_2} &&b_{s_2s_1} \\
b_{s_1s_2s_1} &&b_{s_2s_1s_2} \\
b_{s_1s_2s_1s_2} &&b_{s_2s_1s_2s_1} \\
b_{s_1s_2s_1s_2s_1} &&b_{s_2s_1s_2s_1s_2} \\
\vdots &&\vdots
\end{matrix}
\\ \\
\hbox{$b_w$ basis}
\end{matrix}
$$
$$
\begin{matrix}
\begin{matrix}
&y_{-\alpha_1} \\
y_{\alpha_1} &&y_{-s_2\alpha_1} \\
y_{s_2\alpha_1} &&y_{-s_1s_2\alpha_1} \\
y_{s_1s_2\alpha_1} &&y_{-s_2s_1s_2\alpha_1} \\
y_{s_2s_1s_2\alpha_1} &&y_{-s_1s_2s_1s_2\alpha_1} \\
y_{s_1s_2s_1s_2\alpha_1} &&y_{-s_2s_1s_2s_1s_2\alpha_1} \\
\vdots &&\vdots
\end{matrix}
\\ \\
\hbox{$x_{-\alpha_1}$}
\end{matrix}
\qquad\qquad
\begin{matrix}
\begin{matrix}
&y_{-\alpha_2} \\
y_{-s_1\alpha_2} &&y_{\alpha_2} \\
y_{-s_2s_1\alpha_2} &&y_{s_1\alpha_2} \\
y_{-s_1s_2s_1\alpha_2} &&y_{s_2s_1\alpha_2} \\
y_{-s_2s_1s_2s_1\alpha_2} &&y_{s_1s_2s_1\alpha_2} \\
y_{-s_1s_2s_1s_2s_1\alpha_2} &&y_{s_2s_1s_2s_1\alpha_2} \\
\vdots &&\vdots
\end{matrix}
\\ \\
\hbox{$x_{-\alpha_2}$}
\end{matrix}
$$

Let
\begin{align}
y_{R^-} &= \prod_{\alpha\in R^+} y_{-\alpha}, \\
\Delta_{121}
&= y_{R^-}\left( \frac{1}{y_{-\alpha_2}y_{-\alpha_1}y_{-\alpha_2}}
+\frac{1}{y_{-s_2\alpha_2}y_{-s_2\alpha_1}y_{-\alpha_2}}\right) \nonumber \\
&=
\frac{y_{R^-}}{y_{-\alpha_1}y_{-\alpha_2}y_{-s_1\alpha_2}}
\left( \frac{y_{-s_1\alpha_2}-y_{-\alpha_2}}{y_{-\alpha_1}}+p(y_{\alpha_1},y_{-\alpha_1})y_{-\alpha_2} \right)
 \nonumber \\
\Delta_{212}
&=
\frac{y_{R^-}}{y_{-\alpha_2}y_{-\alpha_1}y_{-s_2\alpha_1}}
\left( \frac{y_{-s_2\alpha_1}-y_{-\alpha_1}}{y_{-\alpha_2}}+p(y_{\alpha_2},y_{-\alpha_2})y_{-\alpha_1}\right),
\quad\hbox{and}
\label{Deltadefn} \\
N&= 1+(1-p(y_{-\alpha_2},y_{-j\alpha_1})y_{-\alpha_2})\big(\sum_{k=1}^{j-1} (1-p(y_{-\alpha_1},y_{-k\alpha_1})y_{-k\alpha_1}\big).
\label{Ndefn}
\end{align}
We note that, \emph{for ordinary cohomology $H_T$ and K-theory $K_T$},
\begin{align*}
N=\begin{cases}
1+(j-1), &\hbox{in $H_T$,} \\
1+e^{-\alpha_2}(e^{-\alpha_1}+\cdots+e^{-(j-1)\alpha_1}), &\hbox{in $K_T$,}
\end{cases}
\qquad\hbox{and}\qquad
\Delta_{121} = \frac{Ny_{R^-}}{y_{-\alpha_1}y_{-\alpha_2}y_{-s_1\alpha_2}}.
\end{align*}

The Schubert and Bott-Samelson cycles for rank 2 and length $\le 1$ are given
$$
\begin{matrix}
\begin{matrix}
&y_{R^-} \\
0 &&0 \\
0 &&0 \\
0 &&0 \\
0 &&0 \\
\vdots &&\vdots
\end{matrix}
\\ \\
[X_1]=[Z_\pt]
\end{matrix}
\qquad\qquad\qquad\qquad
\begin{matrix}
\begin{matrix}
&\frac{y_{R^-}}{y_{-\alpha_1}} \\
\frac{y_{R^-}}{y_{-\alpha_1}} &&0 \\
0 &&0 \\
0 &&0 \\
0 &&0 \\
\vdots &&\vdots
\end{matrix}
\\ \\
[X_{s_1}]=[Z_1]
\end{matrix}
\qquad\qquad\qquad\qquad
\begin{matrix}
\begin{matrix}
&\frac{y_{R^-}}{y_{-\alpha_2}} \\
0 &&\frac{y_{R^-}}{y_{-\alpha_2}} \\
0 &&0 \\
0 &&0 \\
0 &&0 \\
\vdots &&\vdots
\end{matrix}
\\ \\
[X_{s_2}]=[Z_2]
\end{matrix}
$$
The remaining Schubert and Bott-Samelson cycles for rank 2 and length $\le 3$ are given
in Figure \ref{rank2Schuberts}.

\begin{sidewaysfigure}
$$
\begin{matrix}
\begin{matrix}
&\frac{y_{R^-}}{y_{-\alpha_1}y_{-\alpha_2}} \\
\frac{y_{R^-}}{y_{-\alpha_1}y_{-\alpha_2}} &&\frac{y_{R^-}}{y_{-\alpha_2}y_{-s_2\alpha_1}} \\
\frac{y_{R^-}}{y_{-\alpha_2}y_{-s_2\alpha_1}} &&0 \\
0 &&0 \\
0 &&0 \\
\vdots &&\vdots
\end{matrix}
\\ \\
[X_{s_1s_2}]=[Z_{12}]
\end{matrix}
\qquad\qquad
\begin{matrix}
\begin{matrix}
&\frac{y_{R^-}}{y_{-\alpha_1}y_{-\alpha_2}} \\
\frac{y_{R^-}}{y_{-\alpha_1}y_{-s_1\alpha_2}} &&\frac{y_{R^-}}{y_{-\alpha_1}y_{-\alpha_2}} \\
0 &&\frac{y_{R^-}}{y_{-\alpha_1}y_{-s_1\alpha_2}} \\
0 &&0 \\
0 &&0 \\
\vdots &&\vdots
\end{matrix}
\\ \\
[X_{s_2s_1}]=[Z_{21}]
\end{matrix}
$$
$$
\begin{matrix}
\begin{matrix}
&\frac{N y_{R^-}}{y_{-\alpha_1}y_{-\alpha_2}y_{-s_1\alpha_2}}\\
\frac{N y_{R^-}}{y_{-\alpha_1}y_{-\alpha_2}y_{-s_1\alpha_2}} &&\frac{y_{R^-}}{y_{-\alpha_1}y_{-\alpha_2}y_{-s_2\alpha_1}} \\
\frac{y_{R^-}}{y_{-\alpha_1}y_{-\alpha_2}y_{-s_2\alpha_1}} &&\frac{y_{R^-}}{y_{-\alpha_1}y_{-s_1\alpha_2}y_{-s_1s_2\alpha_1}} \\
\frac{y_{R^-}}{y_{-\alpha_1}y_{-s_1\alpha_2}y_{-s_1s_2\alpha_1}} &&0 \\
0 &&0 \\
\vdots &&\vdots
\end{matrix}
\\ \\
[X_{s_1s_2s_1}]
\end{matrix}
\qquad
\begin{matrix}
\begin{matrix}
&\frac{y_{R^-}}{y_{-\alpha_1}y_{-\alpha_2}y_{-s_2\alpha_1}} \\
\frac{y_{R^-}}{y_{-\alpha_1}y_{-\alpha_2}y_{-s_1\alpha_2}} &&\frac{y_{R^-}}{y_{-\alpha_1}y_{-\alpha_2}y_{-s_2\alpha_1}} \\
\frac{y_{R^-}}{y_{-\alpha_2}y_{-s_2\alpha_1}y_{-s_2s_1\alpha_2}} &&\frac{y_{R^-}}{y_{-\alpha_1}y_{-\alpha_2}y_{-s_1\alpha_2}} \\
0 &&\frac{y_{R^-}}{y_{-\alpha_2}y_{-s_2\alpha_1}y_{-s_2s_1\alpha_2}} \\
0 &&0 \\
\vdots &&\vdots
\end{matrix}
\\ \\
[X_{s_2s_1s_2}]
\end{matrix}
$$
$$
\begin{matrix}
\begin{matrix}
&\Delta_{121} \\
\Delta_{121} &&\frac{y_{R^-}}{y_{-\alpha_1}y_{-\alpha_2}y_{-s_2\alpha_1}}  \\
\frac{y_{R^-}}{y_{-\alpha_1}y_{-\alpha_2}y_{-s_2\alpha_1}}
&&\frac{y_{R^-}}{y_{-\alpha_1}y_{-s_1\alpha_2}y_{-s_1s_2\alpha_1}}  \\
\frac{y_{R^-}}{y_{-\alpha_1}y_{-s_1\alpha_2}y_{-s_1s_2\alpha_1}} &&0 \\
0 &&0 \\
\vdots &&\vdots
\end{matrix}
\\ \\
[Z_{121}]
\end{matrix}
\qquad
\begin{matrix}
\begin{matrix}
&\Delta_{212} \\
\frac{y_{R^-}}{y_{-\alpha_1}y_{-\alpha_2}y_{-s_1\alpha_2}} &&\Delta_{212} \\
\frac{y_{R^-}}{y_{-\alpha_2}y_{-s_2\alpha_1}y_{-s_2s_1\alpha_2}} &&\frac{y_{R^-}}{y_{-\alpha_1}y_{-\alpha_2}y_{-s_1\alpha_2}} \\
0 &&\frac{y_{R^-}}{y_{-\alpha_2}y_{-s_2\alpha_1}y_{-s_2s_1\alpha_2}} \\
0 &&0 \\
\vdots &&\vdots
\end{matrix}
\\ \\
[Z_{212}]
\end{matrix}
$$
\caption{Schubert and Bott-Samelson cycles for rank 2 and length $\le 3$.}
\label{rank2Schuberts}
\end{sidewaysfigure}

\subsection{Schubert products in rank 2}

Using the explicit moment graph representations of the Schubert classes, the formulas
for products $g[X_w]$ given at the end of Section \ref{Schubertproducts}
allow for quick computations of the products of Schubert classes in rank 2 for
Weyl group elements up to length 3.
It is straightforward to check that these generalise the corresponding computations
for equivariant cohomology and equivariant K-theory which were given in \cite[\S 5]{GR}.
Since $[X_{s_1s_2s_1s_2}]=[X_{s_2s_1s_2s_1}]=1$ in Type $B_2$, these calculations completely determine
all Schubert products generalized equivariant Schubert products for Types $A_2$ and $B_2$.

The Schubert products  for low dimensional Schubert varieties are as follows.
$$
[X_1]^2 = y_{R^-}[X_1], \qquad
[X_1][X_{s_1}] = \frac{y_{R^-}}{y_{-\alpha_1}}[X_1], \qquad
[X_1][X_{s_2}] = \frac{y_{R^-}}{y_{-\alpha_2}}[X_1], 
$$
$$
[X_1][X_{s_1s_2}] = \frac{y_{R^-}}{y_{-\alpha_1}y_{-\alpha_2}}[X_1], \qquad
[X_1][X_{s_2s_1}] = \frac{y_{R^-}}{y_{-\alpha_2}y_{-\alpha_1}}[X_1], 
$$
$$
[X_1][X_{s_1s_2s_1}] = \frac{N y_{R^-}}{y_{-\alpha_1}y_{-\alpha_2}y_{-s_1\alpha_2}}[X_1], \qquad
[X_1][X_{s_2s_1s_2}] = \frac{y_{R^-}}{y_{-\alpha_2}y_{-\alpha_1}y_{-s_2\alpha_1}}[X_1], 
$$
$$
[X_{s_1}]^2 = \frac{y_{R^-}}{y_{-\alpha_1}}[X_{s_1}], \qquad
[X_{s_1}][X_{s_1s_2}] = \frac{y_{R^-}}{y_{-\alpha_1}y_{-\alpha_2}}[X_{s_1}], \qquad
[X_{s_1}][X_{s_1s_2s_1}] = \frac{N y_{R^-}}{y_{-\alpha_1}y_{-\alpha_2}y_{-s_1\alpha_2}}[X_{s_1}],
$$
\begin{align*}
[X_{s_1}][X_{s_2}] &= \frac{y_{R^-}}{y_{-\alpha_1}y_{-\alpha_2}}[X_1], \\
[X_{s_1}][X_{s_2s_1}] 
&=
\frac{y_{R^-}}{y_{-\alpha_1}y_{-s_1\alpha_2}}[X_{s_1}]+
\frac{y_{R^-}}{y_{-\alpha_2}y_{-\alpha_1}y_{-s_1\alpha_2}}
\left(\frac{y_{-s_1\alpha_2}-y_{-\alpha_2}}{y_{-\alpha_1}}\right) [X_1], \\
[X_{s_1}][X_{s_2s_1s_2}] 
&=\frac{y_{R^-}}{y_{-\alpha_2}y_{-\alpha_1}y_{-s_1\alpha_2}}[X_{s_1}]
+\frac{y_{R^-}}{y_{-\alpha_1}y_{-\alpha_2}y_{-s_1\alpha_2}y_{-s_2\alpha_1}   }
\left(\frac{y_{-s_1\alpha_2}-y_{-s_2\alpha_1}}{y_{-\alpha_1}}\right)[X_1], \\
\end{align*}
$$
[X_{s_2}]^2 = \frac{y_{R^-}}{y_{-\alpha_2}}[X_{s_2}], \quad
[X_{s_2}][X_{s_2s_1}] = \frac{y_{R^-}}{y_{-\alpha_2}y_{-\alpha_1}}[X_{s_2}], \quad
[X_{s_2}][X_{s_2s_1s_2}] = \frac{y_{R^-}}{y_{-\alpha_2}y_{-\alpha_1}y_{-s_2\alpha_1}}[X_{s_2}], 
$$
\begin{align*}
[X_{s_2}][X_{s_1s_2}] &= \frac{y_{R^-}}{y_{-\alpha_2}y_{-s_2\alpha_1}}[X_{s_2}]+
\frac{y_{R^-}}{y_{-\alpha_1}y_{-\alpha_2}y_{-s_2\alpha_1}}
\left(\frac{y_{-s_2\alpha_1}-y_{-\alpha_1}}{y_{-\alpha_2}}\right) [X_1], \\
[X_{s_2}][X_{s_1s_2s_1}] 
&= \frac{y_{R^-}}{y_{-\alpha_1}y_{-\alpha_2}y_{-s_2\alpha_1}}[X_{s_2}]
+\frac{y_{R^-}}{y_{-\alpha_1}y_{-\alpha_2}y_{-s_1\alpha_2}y_{-s_2\alpha_1}}
\left(\frac{Ny_{-s_2\alpha_1}-y_{-s_1\alpha_2}}{y_{-\alpha_2}}\right) [X_1],  
\end{align*}
\begin{align*}
[X_{s_1s_2}]^2
&=\frac{y_{R^-}}{y_{-\alpha_2}y_{-s_2\alpha_1}}[X_{s_1s_2}]
+\frac{y_{R^-}}{y_{-\alpha_2}y_{-\alpha_1}y_{-s_2\alpha_1}}
\left(\frac{y_{-s_2\alpha_1}-y_{-\alpha_1}}{y_{-\alpha_2}}\right) [X_{s_1}],
\\
[X_{s_1s_2}][X_{s_2s_1}]
&= \frac{y_{R^-}}{y_{-\alpha_1}y_{-\alpha_2}y_{-s_1\alpha_2}}[X_{s_1}]
+\frac{y_{R^-}}{y_{-\alpha_1}y_{-\alpha_2}y_{-s_2\alpha_1}}[X_{s_2}] \\
&\qquad+\frac{y_{R^-}}{y_{-\alpha_1}y_{-\alpha_2}y_{-s_1\alpha_2}y_{-s_2\alpha_1}}
\left(\Big(\frac{y_{-s_2\alpha_1}-y_{-\alpha_1}}{y_{-\alpha_2}}\Big)
\Big(\frac{y_{-s_1\alpha_2}-y_{-\alpha_2}}{y_{-\alpha_1}}\Big)-1\right)[X_1], 
\\
[X_{s_1s_2}][X_{s_1s_2s_1}] 
&= \frac{y_{R^-}}{y_{-\alpha_1}y_{-\alpha_2}y_{-s_2\alpha_1}}[X_{s_1s_2}]
+\frac{y_{R^-}}{y_{-\alpha_1}y_{-\alpha_2}y_{-s_1\alpha_2}y_{-s_2\alpha_1}}
\left(\frac{Ny_{-s_2\alpha_1}-y_{-s_1\alpha_2}}{y_{-\alpha_2}}\right) [X_{s_1}], \\  
[X_{s_1s_2}][X_{s_2s_1s_2}] &= 
\frac{y_{R^-}}{y_{-\alpha_2}y_{-s_2\alpha_1}y_{-s_2s_1\alpha_2}}[X_{s_1s_2}] \\
&\qquad+
\frac{y_{R^-}}{y_{-\alpha_1}y_{-\alpha_2}y_{-s_2\alpha_1}y_{-s_1\alpha_2}y_{-s_2s_1\alpha_2}}
\left(\frac{y_{-s_2\alpha_1}y_{-s_2s_1\alpha_2}-y_{-\alpha_1}y_{-s_1\alpha_2}}{y_{-\alpha_2}}
\right)[X_{s_1}] \\
&\qquad+
\frac{y_{R^-}}{y_{-\alpha_1}y_{-\alpha_2}y_{-s_2\alpha_1}y_{-s_2s_1\alpha_2}}
\left(\frac{y_{-s_2s_1\alpha_2}-y_{-\alpha_1}}{y_{-s_2\alpha_1}}\right)[X_{s_2}] \\
&\qquad+
\frac{y_{R^-}}{y_{-\alpha_2}^2}\left(
\frac{1}{y_{-\alpha_1}^2y_{-s_2\alpha_1}}-\frac{1}{y_{-s_2\alpha_1}^2y_{-\alpha_1}}
-\frac{1}{y_{-\alpha_1}^2y_{-s_1\alpha_2}}+\frac{1}{y_{-s_2\alpha_1}^2y_{-s_2s_1\alpha_2}}\right)[X_1], 
\end{align*}
\begin{align*}
[X_{s_2s_1}]^2 
&=\frac{y_{R^-}}{y_{-\alpha_1}y_{-s_1\alpha_2}}[X_{s_2s_1}]
+\frac{y_{R^-}}{y_{-\alpha_1}y_{-\alpha_2}y_{-s_1\alpha_2}}
\left(\frac{y_{-s_1\alpha_2}-y_{-\alpha_2}}{y_{-\alpha_1}}\right) [X_{s_2}], \\
[X_{s_2s_1}][X_{s_1s_2s_1}] &=
\frac{y_{R^-}}{y_{-\alpha_1}y_{-s_1\alpha_2}y_{-s_1s_2\alpha_1}}[X_{s_2s_1}]
+\frac{y_{R^-}}{y_{-\alpha_1}y_{-s_1\alpha_2}^2}\left(\frac{N}{y_{-\alpha_2}}-\frac{1}{y_{-s_1s_2\alpha_1}}\right)[X_{s_1}] \\
&\qquad+\frac{y_{R^-}}{y_{-\alpha_1}^2}\left(
\frac{1}{y_{-s_2\alpha_1}y_{-\alpha_2}}-\frac{1}{y_{-s_1s_2\alpha_1}y_{-s_1\alpha_2}}\right)[X_{s_2}] \\
&\qquad+ \frac{y_{R^-}}{y_{-\alpha_1}^2}\left(
\frac{N}{y_{-\alpha_2}^2y_{-s_1\alpha_2}}
-\frac{N}{y_{-\alpha_2}y_{-s_1\alpha_2}^2}
-\frac{1}{y_{-\alpha_2}^2y_{-s_2\alpha_1}}
+\frac{1}{y_{-s_1\alpha_2}^2y_{-s_1s_2\alpha_1}}
\right)[X_1], \\
[X_{s_2s_1}][X_{s_2s_1s_2}] &= 
\frac{y_{R^-}}{y_{-\alpha_2}y_{-\alpha_1}y_{-s_1\alpha_2}}[X_{s_2s_1}]
+\frac{y_{R^-}}{y_{-\alpha_2}y_{-\alpha_1}^2}\left(\frac{1}{y_{-s_2\alpha_1}}-\frac{1}{y_{-s_1\alpha_2}}\right) [X_{s_2}],
\end{align*}
%
\begin{align*}
&[X_{s_1s_2s_1}]^2 =
\frac{y_{R^-}}{y_{-\alpha_1}y_{-s_1\alpha_2}y_{-s_1s_2\alpha_1}}[X_{s_1s_2s_1}]
+\frac{y_{R^-}}{y_{-\alpha_1^2}}\left(
\frac{1}{y_{-\alpha_2}y_{-s_2\alpha_1}}-\frac{1}{y_{-s_1\alpha_2}y_{-s_1s_2\alpha_1}}\right)[X_{s_1s_2}]
\\
&\quad+\frac{y_{R^-}}{y_{-\alpha_1}y_{-\alpha_2}}\left(
\frac{N^2}{y_{-\alpha_2}y_{-s_1\alpha_2}^2}-\frac{N}{y_{-s_1\alpha_2}^2 y_{-s_1s_2\alpha_1}}
-\frac{1}{y_{-\alpha_1}y_{-\alpha_2}y_{-s_2\alpha_1}}
+\frac{1}{y_{-\alpha_1}y_{-s_1\alpha_2}y_{-s_1s_2\alpha_1}}\right)[X_1],
\\
&[X_{s_1s_2s_1}][X_{s_2s_1s_2}] 
= \frac{y_{R^-}}{y_{-\alpha_1}y_{-\alpha_2}y_{-s_2\alpha_1}y_{-s_2s_1\alpha_2}}[X_{s_1s_2}]
+ \frac{y_{R^-}}{y_{-\alpha_1}y_{-\alpha_2}y_{-s_1\alpha_2}y_{-s_1s_2\alpha_1}}[X_{s_2s_1}] \\
&\qquad\qquad\qquad+ \frac{y_{R^-}}{y_{-\alpha_1}y_{-\alpha_2}}\left(
\frac{N}{y_{-\alpha_2}y_{-s_1\alpha_2}^2}
-\frac{1}{y_{-\alpha_2}y_{-s_2\alpha_1}y_{-s_2s_1\alpha_2}}
-\frac{1}{y_{-s_1\alpha_2}^2y_{-s_1s_2\alpha_1}}\right)[X_{s_1}] \\
&\qquad\qquad\qquad+ \frac{y_{R^-}}{y_{-\alpha_1}y_{-\alpha_2}}\left(
\frac{1}{y_{-\alpha_1}y_{-s_2\alpha_1}^2}
-\frac{1}{y_{-\alpha_1}y_{-s_1\alpha_2}y_{-s_1s_2\alpha_1}}
-\frac{1}{y_{-s_2\alpha_1}^2y_{-s_2s_1\alpha_2}}\right)[X_{s_2}],
\\
&[X_{s_2s_1s_2}]^2 =
\frac{y_{R^-}}{y_{-\alpha_2}y_{-s_2\alpha_1}y_{-s_2s_1\alpha_2}}[X_{s_2s_1s_2}]
+\frac{y_{R^-}}{y_{-\alpha_2^2}}\left(
\frac{1}{y_{-\alpha_1}y_{-s_1\alpha_2}}-\frac{1}{y_{-s_2\alpha_1}y_{-s_2s_1\alpha_2}}\right)[X_{s_2s_1}]
\\
&\quad+\frac{y_{R^-}}{y_{-\alpha_1}y_{-\alpha_2}}\left(
\frac{1}{y_{-\alpha_1}y_{-s_2\alpha_1}^2}-\frac{1}{y_{-s_2\alpha_1}^2 y_{-s_2s_1\alpha_2}}
-\frac{1}{y_{-\alpha_1}y_{-\alpha_2}y_{-s_1\alpha_2}}
+\frac{1}{y_{-\alpha_2}y_{-s_2\alpha_1}y_{-s_2s_1\alpha_2}}\right)[X_1],
\end{align*}

\section{The calculus of BGG operators}\label{BGGopcalc}

The \emph{nil affine Hecke algebra} is the algebra over $\LL$ with generators $x_\lambda$, $y_\lambda$,
$t_w$, with $\lambda,\mu\in \fh_\ZZ^*$ and $w\in W_0$, with relations
$$x_{\lambda+\mu} = x_\lambda+x_\mu-p(x_\lambda,x_\mu)x_\lambda x_\mu,
\qquad
y_{\lambda+\mu} = y_\lambda+y_\mu-p(y_\lambda,y_\mu)y_\lambda y_\mu,
\qquad
x_\lambda y_\mu = y_\mu x_\lambda, 
$$
and
$$t_vt_w = t_{vw}, \qquad
t_w y_\lambda = y_\lambda t_w, \qquad
t_w x_\lambda = x_{w\lambda} t_w,
\qquad\hbox{for $v, w\in W_0$, $\lambda\in \fh_\ZZ^*$.}
$$

Recall from \eqref{pushpull} that the \emph{pushpull operators}, or \emph{BGG-Demazure operators} 
are given by
\begin{equation}
A_i = (1+t_{s_i})\frac{1}{x_{-\alpha_i}},
\qquad\hbox{for $i=1,2,\ldots, n$.}
\end{equation}
In general,
\begin{align}
A_i &= (1+t_{s_i})\frac{1}{x_{-\alpha_i}} = \frac{1}{x_{-\alpha_i}}+\frac{1}{x_{\alpha_i}}t_{s_i}
=\frac{1}{x_{-\alpha_i}}-\frac{1-p(x_{\alpha_i},x_{-\alpha_i})x_{-\alpha_i}}{x_{-\alpha_i}}t_{s_i}  \nonumber \\
&=\frac{1}{x_{-\alpha_i}}\left(1-(1-p(x_{\alpha_i},x_{-\alpha_i})x_{-\alpha_i})t_{s_i}\right)
=\frac{1}{x_{-\alpha_i}}(1-t_{s_i})+p(x_{\alpha_i},x_{-\alpha_i})t_{s_i}.
\end{align}
so that $A_i$ is a divided difference operator plus an extra term.
As in \cite[Prop.\ 3.1]{BE1},
\begin{align*}
A_i^2 &= (1+ t_{s_i})\frac{1}{x_{-\alpha_i}}(1+ t_{s_i})\frac{1}{x_{-\alpha_i}}
=\left(\frac{1}{x_{-\alpha_i}}+\frac{1}{x_{\alpha_i}}t_{s_i}\right)
(1+t_{s_i})\frac{1}{x_{-\alpha_i}} \\
&=\left(\frac{1}{x_{-\alpha_i}}+\frac{1}{x_{\alpha_i}}\right)
(1+t_{s_i})\frac{1}{x_{-\alpha_i}}
=\left(\frac{1}{x_{-\alpha_i}}+\frac{1}{x_{\alpha_i}}\right)
A_i,
\end{align*}
so that
\begin{equation}
A_i^2 = 
\left(\frac{1}{x_{-\alpha_i}}+\frac{1}{x_{\alpha_i}}\right)
A_i  = 
A_i
\left(\frac{1}{x_{-\alpha_i}}+\frac{1}{x_{\alpha_i}}\right)
=A_ip(x_{\alpha_i},x_{-\alpha_i}).
\end{equation}
Note also that
\begin{align}
t_{s_i}A_i &= t_{s_i}(1+t_{s_i})\frac{1}{x_{-\alpha_i}}=A_i
\quad\hbox{and}\quad \\
A_it_{s_i} &= (1+t_{s_i})\frac{1}{x_{-\alpha_i}}t_{s_i}
= (1+t_{s_i})\frac{1}{x_{\alpha_i}}
=A_i\frac{x_{-\alpha_i}}{x_{\alpha_i}}.
\end{align}
If $f\in \LL[[x_\lambda\ |\ \lambda\in \fh_\ZZ^*]]$ then
\begin{align*}
fA_i &= f(1+t_{s_i})\frac{1}{x_{-\alpha_i}}
=f \frac{1}{x_{-\alpha_i}}+ft_{s_i}\frac{1}{x_{-\alpha_i}}
\quad\hbox{and}\quad \\
A_i(s_if) &= (1+t_{s_i})\frac{s_if}{x_{-\alpha_i}}
=(s_if+ft_{s_i})\frac{1}{x_{-\alpha_i}},
\end{align*}
so that
\begin{equation}\label{affHeckerelation}
fA_i = A_i(s_i f)+ \left(\frac{f-s_if}{x_{-\alpha_i}}\right).
\end{equation}
The relation \eqref{affHeckerelation} is the analogue, for this setting, of a key relation in the
definition of the classical nil-affine Hecke algebra (see \cite[Lemma 7.1.10]{CG} or \cite[(1.3)]{GR}).

Next are useful, expansions of products of $t_{s_i}$ in terms of products of $A_i$ with $x$s on the left,
\begin{align*}
t_{s_1} &= x_{\alpha_1}A_1 - \frac{x_{\alpha_1}}{x_{-\alpha_1}}, \\
t_{s_2}t_{s_1} 
&=x_{s_2\alpha_1}
x_{\alpha_2}A_2A_1 - x_{s_2\alpha_1}\frac{x_{\alpha_2}}{x_{-\alpha_2}}A_1
-\frac{x_{s_2\alpha_1}}{x_{-s_2\alpha_1}}x_{\alpha_2}A_2 
+ \frac{x_{s_2\alpha_1}}{x_{-s_2\alpha_1}}\frac{x_{\alpha_2}}{x_{-\alpha_2}}\\
t_{s_1}t_{s_2}t_{s_1}
&=
x_{s_1s_2\alpha_1}x_{s_1\alpha_2}x_{\alpha_1}A_1A_2A_1
-x_{s_1s_2\alpha_1}x_{s_1\alpha_2}\frac{x_{\alpha_1}}{x_{-\alpha_1}}A_2A_1
-\frac{x_{s_1s_2\alpha_1}}{x_{-s_1s_2\alpha_1}}x_{s_1\alpha_2}x_{\alpha_1}A_1A_2 \\
&\qquad +\frac{x_{s_2s_1\alpha_2}}{x_{-s_2s_1\alpha_2}}x_{s_2\alpha_1}
\frac{x_{\alpha_2}}{x_{-\alpha_2}}A_1
+\frac{x_{s_1s_2\alpha_1}}{x_{-s_1s_2\alpha_1}}x_{s_1\alpha_2}
\frac{x_{\alpha_1}}{x_{-\alpha_1}}A_2 
- \frac{x_{s_1s_2\alpha_1}}{x_{-s_1s_2\alpha_1}}\frac{x_{s_1\alpha_2}}{x_{-s_1\alpha_2}}\frac{x_{\alpha_1}}{x_{-\alpha_1}} \\
&\qquad
+\left(\frac{x_{s_1\alpha_2}}{x_{-s_1\alpha_2}}
\frac{x_{s_1s_2\alpha_1}}{x_{-s_1s_2\alpha_1}}x_{\alpha_1}
- \frac{x_{s_1\alpha_2}}{x_{-s_1\alpha_2}}x_{s_1s_2\alpha_1} 
-\frac{x_{s_2s_1\alpha_2}}{x_{-s_2s_1\alpha_2}}x_{s_2\alpha_1}
\frac{x_{\alpha_2}}{x_{-\alpha_2}}
\right)A_1 \\
t_{s_1}t_{s_2}t_{s_1}t_{s_2}
&=
x_{s_2s_1s_2\alpha_1}x_{s_2s_1\alpha_2}x_{s_2\alpha_1}x_{\alpha_2}A_2A_1A_2A_1 \\
&\qquad
-x_{s_2s_1s_2\alpha_1}x_{s_2s_1\alpha_2}x_{s_2\alpha_1}\frac{x_{\alpha_2}}{x_{-\alpha_2}}A_1A_2A_1
-\frac{x_{s_2s_1s_2\alpha_1}}{x_{-s_2s_1s_2\alpha_1}}x_{s_2s_1\alpha_2}x_{s_2\alpha_1}x_{\alpha_2}A_2A_1A_2 \\
&\qquad
+\frac{x_{s_2s_1s_2\alpha_1}}{x_{-s_2s_1s_2\alpha_1}}x_{s_2s_1\alpha_2}x_{s_2\alpha_1}\frac{x_{\alpha_2}}{x_{-\alpha_2}}A_1A_2 \\
&\qquad
+\left(
\frac{x_{s_2s_1s_2\alpha_1}}{x_{-s_2s_1s_2\alpha_1}}\frac{x_{s_2s_1\alpha_2}}{x_{-s_2s_1\alpha_2}}x_{s_2\alpha_1}x_{\alpha_2}
- x_{s_2s_1s_2\alpha_1}\frac{x_{s_2s_1\alpha_2}}{x_{-s_2s_1\alpha_2}}x_{\alpha_2}
-x_{s_2s_1s_2\alpha_1}x_{s_2s_1\alpha_2}\frac{x_{s_2\alpha_1}}{x_{-s_2\alpha_1}}
\right)A_2A_1 \\
&\qquad
-\left(
\frac{x_{s_2s_1s_2\alpha_1}}{x_{-s_2s_1s_2\alpha_1}}\frac{x_{s_2s_1\alpha_2}}{x_{-s_2s_1\alpha_2}}x_{s_2\alpha_1}
- x_{s_2s_1s_2\alpha_1}\frac{x_{s_2s_1\alpha_2}}{x_{-s_2s_1\alpha_2}}\right)\frac{x_{\alpha_2}}{x_{-\alpha_2}}A_1 \\
&\qquad
+\left(
\frac{x_{s_2s_1s_2\alpha_1}}{x_{-s_2s_1s_2\alpha_1}}x_{s_2s_1\alpha_2}
\frac{x_{s_2\alpha_1}}{x_{-s_2\alpha_1}}
- \frac{x_{s_2s_1s_2\alpha_1}}{x_{-s_2s_1s_2\alpha_1}}\frac{x_{s_2s_1\alpha_2}}{x_{-s_2s_1\alpha_2}}\frac{x_{s_2\alpha_1}}{x_{-s_2\alpha_1}}x_{\alpha_2}\right)A_2 \\
&\qquad
+ \frac{x_{s_2s_1s_2\alpha_1}}{x_{-s_2s_1s_2\alpha_1}}\frac{x_{s_2s_1\alpha_2}}{x_{-s_2s_1\alpha_2}}\frac{x_{s_2\alpha_1}}{x_{-s_2\alpha_1}}\frac{x_{\alpha_2}}{x_{-\alpha_2}},
\end{align*}
and expansions of products of $t_{s_i}$ in terms of products of $A_i$ with $x$s on the right,
\begin{align*}
t_{s_1} 
&= A_1x_{-\alpha_1}-1, \\
t_{s_1}t_{s_2} 
&=A_1A_2x_{-\alpha_2}x_{-s_2\alpha_1}-A_1x_{-s_2\alpha_1}-A_2x_{-\alpha_2}+1, \\
t_{s_1}t_{s_2}t_{s_1}
&=A_1A_2A_1x_{-\alpha_1}x_{-s_1\alpha_2}x_{-s_1s_2\alpha_1}-A_1A_2x_{-s_1\alpha_2}x_{-s_1s_2\alpha_1}
-A_2A_1x_{-\alpha_1}x_{-s_1\alpha_2} \\
&\qquad
+A_1x_{-s_2\alpha_1}+A_2x_{-s_1\alpha_2}-1 
+A_1\left(x_{-\alpha_1} -x_{-s_2\alpha_1}
-\frac{x_{-\alpha_1}}{x_{\alpha_1}}x_{-s_1s_2\alpha_1}\right), \\
t_{s_1}t_{s_2}t_{s_1}t_{s_2}
&=A_1A_2A_1A_2x_{-\alpha_2}x_{-s_2\alpha_1}x_{-s_2s_1\alpha_2}x_{-s_2s_1s_2\alpha_1} \\
&\qquad-A_1A_2A_1x_{-s_2\alpha_1}x_{-s_2s_1\alpha_2}x_{-s_2s_1s_2\alpha_1}
-A_2A_1A_2x_{-\alpha_2}x_{-s_2\alpha_1}x_{-s_2s_1\alpha_2} \\
&\qquad+A_1A_2\left(
-\frac{x_{-\alpha_2}}{x_{\alpha_2}}x_{-s_2s_1\alpha_2}x_{-s_2s_1s_2\alpha_1}
-x_{-\alpha_2}\frac{x_{-s_2\alpha_1}}{x_{s_2\alpha_1}}x_{-s_2s_1s_2\alpha_1}
+x_{-\alpha_2}x_{-s_2\alpha_1}\right) \\
&\qquad+A_2A_1x_{-s_2\alpha_1}x_{-s_2s_1\alpha_2} \\
&\qquad-A_1\left(x_{-s_2\alpha_1}
-\frac{x_{-s_2\alpha_1}}{x_{s_2\alpha_1}}x_{-s_2s_1s_2\alpha_1}\right)
-A_2\left(x_{-\alpha_2} - \frac{x_{-\alpha_2}}{x_{\alpha_2}}x_{-s_2s_1\alpha_2}\right) + 1.
\end{align*}
Finally, there are expansions of products of $A_i$ in terms of products of $t_{s_i}$:
\begin{align*}
A_1 &= (t_{s_1}+1)\frac{1}{x_{-\alpha_1}}, \\
A_1A_2 
&= (t_{s_1}+1)\left(t_{s_2}\frac{1}{x_{-\alpha_2}x_{-s_2\alpha_1}}+\frac{1}{x_{-\alpha_1}x_{-\alpha_2}}\right), \\
A_1A_2A_1
&=(t_{s_1}+1)\left(
\begin{matrix}
\displaystyle{
t_{s_2}t_{s_1}\frac{1}{x_{-\alpha_1}x_{-s_1\alpha_2}x_{-s_1s_2\alpha_1}}
+t_{s_2}\frac{1}{x_{-\alpha_1}x_{-\alpha_2}x_{-s_2\alpha_1}} } \\ \\
\displaystyle{
\qquad+\frac{1}{x_{-\alpha_1}}\bigg(\frac{1}{x_{-\alpha_1}x_{-\alpha_2}}+\frac{1}{x_{-s_1\alpha_1}x_{-s_2\alpha_1}}\bigg) }
\end{matrix}
\right),\\
A_1A_2A_1A_2
&=(t_{s_1}+1)\left(
\begin{matrix}
\displaystyle{
t_{s_2}t_{s_1}t_{s_2}\frac{1}{x_{-\alpha_2}x_{-s_2\alpha_1}x_{-s_2s_1\alpha_2}x_{-s_2s_1s_2\alpha_1}} } \hfill \\ \\
\displaystyle{
+t_{s_2}t_{s_1}\frac{1}{x_{-\alpha_2}x_{-\alpha_1}x_{-s_1\alpha_2}x_{-s_1s_2\alpha_1}} } \hfill \\ \\
\displaystyle{
+t_{s_2}\frac{1}{x_{-\alpha_2}x_{-s_2\alpha_1}}
\bigg(\frac{1}{x_{-\alpha_2}x_{-\alpha_1}}+\frac{1}{x_{-s_2\alpha_1}x_{-s_2\alpha_2}}
+\frac{1}{x_{-s_2s_1\alpha_2}x_{-s_2s_1\alpha_1}}\bigg) } \\ \\
\displaystyle{
+\frac{1}{x_{-\alpha_1}x_{-\alpha_2}}\bigg(
\frac{1}{x_{-\alpha_2}x_{-\alpha_1}}+\frac{1}{x_{-s_2\alpha_1}x_{-s_2\alpha_2}}+\frac{1}{x_{-s_1\alpha_2}x_{-s_1\alpha_1}}
\bigg) } \hfill
\end{matrix}
\right).
\end{align*}
These formulas arranged so that products beginning with $t_{s_2}$ and $A_2$ are obtained from the above
formulas by switching 1s and 2s.  In particular, the ``braid relations'' for the operators $A_i$ are the 
equations given by, for example, in the case that $s_1s_2s_1=s_2s_1s_2$ so that
$s_1\alpha_2 = s_2\alpha_1 = \alpha_1+\alpha_2$ then 
$$0 = t_{s_1}t_{s_2}t_{s_1} - t_{s_2}t_{s_1}t_{s_2}$$
is equivalent to 
\begin{align*}
A_2A_1A_2 &- \left(\frac{1}{x_{-\alpha_2}x_{-\alpha_1}}
-\frac{1}{x_{-\alpha_1}x_{-\alpha_3}}
+\frac{1}{x_{\alpha_2}x_{-\alpha_3}}\right) A_2 \\
&=A_1A_2A_1 - \left(\frac{1}{x_{-\alpha_1}x_{-\alpha_2}}
-\frac{1}{x_{-\alpha_2}x_{-\alpha_3}}
+\frac{1}{x_{\alpha_1}x_{-\alpha_3}}\right) A_1,
\end{align*}
as indicated in \cite[Proposition 5.7]{SZ}.


\begin{thebibliography}{99}

\bibitem[AB]{AB}  M.F.\ Atiyah and R.\  Bott, \emph{The moment map and equivariant cohomology}, Topology \textbf{23} (1984)1--28, MR 0721448.

\bibitem[Ad]{Ad} J.F. Adams, \emph{Stable homotopy and generalized homology},
Univ. of Chicago Press, 1974, MR1324104.

\bibitem[An]{An} M.\ Ando, \emph{Power operations in elliptic cohomology and representations of loop groups}, 
Trans.\ Amer.\ Math.\ Soc.\ \textbf{352} (2000) 5619Ð5666, MR1637129 arXiv:


\bibitem[AGM]{AGM} D.\ Anderson, S.\  Griffeth, and E.\  Miller, \emph{Positivity and Kleiman transversality in equivariant K-theory of homogeneous spaces}, J.\ Eur.\ Math.\ Soc.\  \textbf{13} (2011) 57-84. MR2735076, arXiv:0808.2785

\bibitem[AS]{AS} D.\ Anderson and A.\ Stapledon, \emph{Schubert varieties are lag Fano over the integers}, to appear in 
Proc.\ Amer.\ Math.\ Soc.\ arXiv:1203.6678

\bibitem[BS]{BS} N.\ Bergeron and F.\ Sottile, \emph{Schubert polynomials, the Bruhat order, and the geometry of flag manifolds}, Duke Math.\ J.\ \textbf{95} (1998) 373--423, MR1652021, arXiv:9703001.

%
%

\bibitem[BiLa]{BiLa} S.\ Billey and V. Lakshmibai, \emph{Singular loci of Schubert varieties}, Progress in Mathematics
\textbf{182} Birkh\"auser Boston, 2000. xii+251 pp. ISBN: 0-8176-4092-4, MR1782635.

\bibitem[Bo]{Bo} A.\ Borel, \emph{Sur la cohomologie des espaces fibr\'es principaux et des espaces
homog\`enes de groupes de Lie compacts}, Ann.\ of Math.\ (2) \textbf{57} (1953) 115--207, MR0051508.

\bibitem[BL]{BL} L.\ Borisov and A. Libgober, \emph{Elliptic genera of singular varieties}, Duke Math.\ J.\
\textbf{116} (2003) 319--351, MR1953295, arXiv:0007108.

%

\bibitem[BE1]{BE1} P. Bressler and S. Evens,
\emph{The Schubert calculus, braid relations, and generalized cohomology},
Trans.\ Amer.\ Math.\ Soc.\ \textbf{317} (1990), 799--811, MR0968883

\bibitem[BE2]{BE2} P. Bressler and S. Evens,
\emph{Schubert calculus in complex cobordism},
Trans.\ Amer.\ Math.\ Soc.\ \textbf{331} (1992), 799--813, MR1044959.

\bibitem[Ch]{Ch} W.-L.\ Chow, \emph{Algebraic varieties with rational dissections},
Proc.\ Nat.\ Acad.\ Sci.\ U.S.A.\ \textbf{42} (1956) 116--119, MR0078006.

\bibitem[CG]{CG} N.\ Chriss and V. Ginzburg,  
{\sl Representation theory and complex geometry}, Birkh\"auser Boston, Inc., Boston, MA, 1997. x+495 pp. 
ISBN: 0-8176-3792-3, MR1433132 and MR2838836.

\bibitem[CPZ]{CPZ} B.\ Calm\`es, V.\ Petrov, K. Zainoulline,
\emph{Invariants, torsion indices and oriented cohomology of complete flags},
arxiv:0905.1341


\bibitem[EW]{EW} B.\ Elias and G.\ Williamson, \emph{The Hodge theory for Soergel bimodules}, in preparation, November 2012.

\bibitem[Fu]{Fu} W.\ Fulton, {\sl Intersection Theory}, Second edition, Ergebnisse der Mathematik und ihrer Grenzgebiete 3 Folge, Springer-Verlag, Berlin, 1998. xiv+470 pp. ISBN: 3-540-62046-X; 0-387-98549-2, MR1644323.

\bibitem[Ga]{Ga} N.\ Ganter, \emph{The elliptic Weyl character formula}, arXiv:1206.0528.

\bibitem[GKM]{GKM}  M. Goresky, R.\ Kottwitz, R.\ MacPherson, \emph{ Equivariant cohomology, Koszul duality, and the localization theorem}, Invent.\ Math.\ 131 (1998), no. 1, 25--83, MR1489894.

\bibitem[GKV]{GKV} V. Ginzburg, M. Kapranov and E. Vasserot, \emph{Elliptic algebras and 
equivariant elliptic cohomology I}, arXiv:9505012.

\bibitem[Gra]{Gra} W.\ Graham, \emph{Positivity in equivariant Schubert calculus}, Duke Math.\ J.\ \textbf{109} (2001) 599-614.
MR1853356 arXiv:9908172

\bibitem[GR]{GR} S.\ Griffeth and A.\ Ram, \emph{Affine Hecke algebras and the Schubert calculus},
European J.\ Combinatorics \textbf{25} (2004) 1263--1283, MR2095481, arXiv:0405333

\bibitem[Gr]{Gr} I. Grojnowski, \emph{Delocalised equivariant elliptic cohomology}, preprint 1994, 
appeared in {\sl Elliptic cohomology} 114--121, London Math.\ Soc.\ Lecture Note Ser.\ \textbf{342} Cambridge Univ.\ Press, Cambridge, 2007, available from http://www.dpmms.cam.ac.uk/$\sim$groj/papers.html.

\bibitem[G]{G} A.\ Grothendieck, \emph{Sur quelques propri\'et\'es fondamentales en th\'eorie des intersections},
S\'eminaire Claude Chevalley \textbf{3} (1958), Expos\'e No. 4, 36 p. 

\bibitem[HHH]{HHH} M.\ Harada, A.\ Henriques, T. Holm, \emph{Computation of 
generalized equivariant cohomologies of Kac-Moody flag varieties}, Adv. Math. \textbf{197}
(2005) 198--221, MR2166181, arXiv:0409305

\bibitem[HLSZ]{SZ} 
A.\  Hoffnung, J.\ Lopez, A.\ Savage, K.\ Zainouilline,
\emph{Formal Hecke algebras and algebraic oriented cohomology theories},
arXiv:1208.4114.

\bibitem[HK]{HK} J.\ Hornbostel and V.\ Kiritchenko, \emph{Schubert calculus for
algebraic cobordism},  J.\ Reine Angew.\ Math.\ \textbf{656} (2011) 59--85, MR2818856, arxiv:0903.3936.



\bibitem[Ka]{Ka} S.\ Kaji, \emph{Schubert calculus, seen from torus equivariant topology}, Trends in Mathematics-New Series
\textbf{12} (2010), 71--89.

\bibitem[KiKr]{KiKr} V.\ Kiritchenko and A.\ Krishna, \emph{Equivariant cobordism of flag varieties and of symmetric varieties},
arXiv:1104.1089.

\bibitem[KK1]{KK1} B.\ Kostant and S.\ Kumar, \emph{The nil Hecke ring and cohomology of $G/P$ for a Kac-Moody group $G$},
Adv.\ in Math.\ \textbf{62} (1986) 187-237, MR0866159.

\bibitem[KK2]{KK2} B.\ Kostant and S.\ Kumar, \emph{$T$-equivariant $K$-theory of
generalized flag manifolds}, Adv.\ in Math.\ \textbf{62} (1986) 187--237, MR0895705.

\bibitem[KL]{KL} D.\ Kazhdan and G.\ Lusztig, \emph{Proof of the Deligne-Langlands conjecture for Hecke algebras}, 
Invent.\ Math.\ \textbf{87} (1987) 153--215, MR0862716.

\bibitem[KP]{KP} V.\ Kac and D.\ Peterson, \emph{Infinite dimensional Lie algebras, 
theta functions and modular forms},  Adv.\ in Math.\ \textbf{53} (1984), 125--264, MR0750341.

\bibitem[Kr]{Kr} A.\ Krishna, \emph{Equivariant cobordism of schemes}, Doc.\ Math.\ \textbf{17} (2012) 95--134, MR2889745, arXiv: 1006.3176

\bibitem[Ku]{Ku} S.\ Kumar, \emph{The nil Hecke ring and singularity of Schubert varieties}, Invent.\ Math.\ \textbf{123} (1996) 471--506, MR1383959.

\bibitem[Ku2]{Ku2} S.\ Kumar, {\sl Kac-Moody groups, their flag varieties and representation theory},
Progress in Mathematics \textbf{204} Birkh\"auser Boston, Inc., Boston, MA, 2002 ISBN: 0-8176-4227-7, MR1923198.

\bibitem[KS]{KS} S.\ Kumar and K.\ Schwede, \emph{Richardson varieties have Kawamata log terminal singularities},
arXiv:1203.6126


\bibitem[LM]{LM} M.\ Levine and F.\ Morel, {\sl Algebraic Cobordism}, Springer Monographs in Mathematics, Springer, Berlin, 2007. xii+244 pp. ISBN: 978-3-540-36822-9; 3-540-36822-1,
MR2286826.

\bibitem[LSS]{LSS} T. Lam, A. Schilling and M. Shimozono, 
\emph{K-theory Schubert calculus of the affine Grassmannian}, 
Compos.\ Math.\ \textbf{146} (2010) 811--852. MR2660675



\bibitem[Lu]{Lu} J.\ Lurie, \emph{A survey of elliptic cohomology}, Algebraic topology, 219-277, Abel Symp.\ \textbf{4} Springer, Berlin, 2009. MR2597740, arXiv:


\bibitem[Ma]{Ma} J.P.\ May, \emph{Equivariant homotopy and cohomology theory}, CBMS Regional Conference Series in Mathematics \textbf{91}. Published for the Conference Board of the Mathematical Sciences, Washington, DC; by the American Mathematical Society, Providence, RI, 1996, ISBN: 0-8218-0319-0, MR1413302.

\bibitem[MR]{MR}  {\sl Elliptic cohomology. Geometry, applications, and higher chromatic analogues}, Papers from the Workshop on Elliptic Cohomology and Chromatic Phenomena held in Cambridge, December 9Ð20, 2002. Edited by Haynes R.\ Miller and Douglas C.\ Ravenel, London Mathematical Society Lecture Note Series \textbf{342} Cambridge University Press, Cambridge, 2007. xii+364 pp. ISBN: 978-0-521-70040-5; 0-521-70040-X MR2330502

\bibitem[Oko]{Oko} C.\ Okonek,\emph{Der Conner-Floyd-Isomorphismus fŸr Abelsche Gruppen}, 
Math.\ Z.\ \textbf{179} (1982) 201-212, MR0645496.



\bibitem[Soe]{Soe} W.\ Soergel, \emph{Kazhdan-Lusztig-Polynome und unzerlegbare Bimoduln Ÿber Polynomringen},
J.\ Inst.\ Math.\ Jussieu \textbf{6} (2007) 501Ð525. MR2329762


\bibitem[To]{To} B.\ Totaro, \emph{The elliptic genus of a singular variety}, in {\sl Elliptic cohomology}, 360--364, London Math.\ Soc.\ Lecture Note Ser.\ \textbf{ 342} Cambridge Univ. Press, Cambridge, 2007, MR2330522

\bibitem[Ty]{Ty} J.\ Tymoczko, \emph{Permutation representations on Schubert varieties}, Amer.\ J.\ Math.\ \textbf{130} (2008) 1171--1194, MR2450205, arXiv: 0604578.

%


\end{thebibliography}
\end{document}